\documentclass[a4paper,11pt]{scrartcl}%
\usepackage{mathtools}
\mathtoolsset{showonlyrefs}
\usepackage{amsmath, amsthm, amssymb}
\usepackage{fontenc}
\usepackage[utf8]{inputenc}
\usepackage[english]{babel}
\usepackage{graphicx}
\usepackage{subcaption,xargs}
\usepackage{algorithm, booktabs}
\usepackage[noend]{algpseudocode}
\usepackage[numbers]{natbib}
\usepackage{enumitem,listings}
\usepackage{hyperref}
\usepackage{multirow,environ,xcolor}
\setkomafont{sectioning}{\rmfamily\bfseries}
\setkomafont{title}{\rmfamily}
\newtheorem{theorem}{Theorem}[section]
\newtheorem{definition}[theorem]{Definition}

\newtheorem{lemma}[theorem]{Lemma}

\newtheorem{example}[theorem]{Example}

\DeclareCaptionLabelSeparator{periodspace}{.\ }
\DeclareCaptionLabelSeparator{periodspace}{.\ }
\title{A Framework for FFT-based Homogenization on Anisotropic Lattices}

\author{%
Ronny Bergmann\footnotemark[1]%
\and
Dennis Merkert\footnotemark[1]%
}
\newcommand{\e}{{\,\mathrm{e}}}
\newcommand{\im}{{\mathrm{i}}}
\newcommand{\dx}{\,\mathrm{d}}
\makeatletter
\newcommandx{\abs}[2][1=\@empty]{#1\lvert #2 #1\rvert}
\newcommandx{\norm}[3][1=\@empty,3=\@empty]{#1\lVert #2 #1\rVert_{#3}}
\makeatother
\newcommand{\tT}{\mathrm{T}}
\newcommand{\hH}{\mathrm{H}}
\newcommand{\vect}[1]{\mathbf{#1}}
\newcommand{\mat}[1]{\mathbf{#1}}
\newcommand{\Stiffness}{\mathcal{C}}
\newcommand{\eff}{\mathrm{eff}}

\newcommand{\Z}{\mathbb{Z}}
\newcommand{\R}{\mathbb{R}}

\newcommand{\T}{\mathbb{T}}
\newcommand{\tensorProd}{:}
\DeclareMathOperator{\Div}{div}
\DeclareMathOperator{\Id}{Id}
\DeclareMathOperator{\Fourier}{\mathcal F}

\DeclareMathOperator{\generatingSet}{\mathcal G}
\DeclareMathOperator{\Grad}{\nabla}
\DeclareMathOperator{\Pattern}{\mathcal P}
\DeclareMathOperator{\GradSym}{\Grad_{\mathrm{Sym}}}

\DeclareMathOperator{\diag}{diag}

\hypersetup{pdfauthor={Ronny Bergmann, Dennis Merkert},
	pdftitle={A Framework for FFT-based Homogenization on Anisotropic Lattices},
	pdfsubject={preprint},%
	pdfcreator = {pdflatex and TextMate},%
	pdfkeywords={},
	plainpages=false, pdfstartview=FitH, pdfview=FitH, pdfpagemode=UseOutlines,%
	bookmarksnumbered=true,bookmarksopen=false,bookmarksopenlevel=0,%
	colorlinks=true,linkcolor=blue!50!black,citecolor=green!50!black,urlcolor=red!50!black%
}
\date{June 21, 2016}
\begin{document}
\maketitle
\footnotetext[1]{%
Department of Mathematics,
Technische Universität Kaiserslautern,
Kaiserslautern, Germany,
\\$\{$bergmann, dmerkert$\}$@mathematik.uni-kl.de.}
\begin{abstract}
	\textbf{Abstract.} In order to take structural anisotropies of a given composite and different shapes of its unit cell into account, we generalize the Basic Scheme in homogenization by Moulinec and Suquet 
	to arbitrary sampling lattices and tilings of the \(d\)-dimensional Euclidean space.
	We employ a Fourier transform for these lattices by introducing the corresponding set of sample points, the so called pattern, and its frequency
	set, the generating set, both representing the anisotropy of both the shape of the unit cell and the chosen preferences in certain sampling directions.
	In several cases, this Fourier transform is of lower dimension than the
	space itself. For the so called rank-\(1\)-lattices it even reduces to a
	one-dimensional Fourier transform having the same leading coefficient as the
	fastest Fourier transform implementation available.
	We illustrate the generalized
	Basic Scheme on an anisotropic laminate and a generalized ellipsoidal Hashin
	structure. For both we give an analytical solution to the elasticity problem,
	in two- and three dimensions, respectively.
	We then illustrate the possibilities of choosing a pattern.
	Compared to classical grids this introduces both a reduction of
	computation time and a reduced error of the numerical method. It also allows
	for anisotropic subsampling, i.e.~choosing a sub lattice of a pixel or voxel
	grid based on anisotropy information of the material at hand.
	\\[\baselineskip]
	\\[\baselineskip]
\end{abstract}
%

%
\section{Introduction}\label{sec:Introduction}

Modern materials are often composites of multiple components which are designed
to obtain overall properties like high durability, flexibility or stiffness.
These inhomogeneities are usually small in comparison to the overall structure
of the material or tool.
Therefore it is computationally beneficial and sometimes even necessary to
replace the inhomogeneous material by a homogeneous one having the same
macroscopic properties, called homogenization.

The underlying assumption is that the microstructure can be represented by a
reference volume that can be repeated periodically to generate the geometry.
While many of these microstructures show macroscopically isotropic behavior
there are also composites that have one or multiple predominant directions.

The classical algorithm to solve such homogenization
problems of periodic microstructures on regular grids was proposed by
Moulinec and Suquet~\cite{MoulinecSuquet1994,MoulinecSuquet1998}.
This algorithm is also called the Basic Scheme and has evoked many
enhancements and modifications. Amongst them are different discretization
methods of the differential operator~\cite{Willot2015,Schneider2015,Schneider2016},
adaptions to composites with infinite contrast, e.g.~porous
media,~\cite{Michel2000,Schneider2016}, incorporation of additional information
about the geometry~\cite{KMS:2015Homogenization} and the solution of
homogenization problems of higher order~\cite{Tran2012}.

All of them have in common that they are formulated on regular tensor product
grids, i.e.~they make use of the commonly known multidimensional Fast Fourier
Transform (FFT). In some cases it is not possible to rotate the
representative volume element without violating the periodicity condition of
the microstructure and then a change of the discretization grid can remedy this.

Galipeau and Casta\~neda~\cite{Galipeau2013,Galipeau2013b}, for example,
construct a periodic laminate structure of elastomers where each of the two
phases consists of aligned elongated particles of a magnetic material.
In this material, the two phases differ in orientation and do not face into the
direction of lamination nor orthogonal to it.
Lahellec, et.~al.~\cite{Lahellec2003} consider a multi-particle problem where
they have an evolving computational grid. The basis vectors of this grid depend
on the macroscopic velocity of a Newtonian fluid and they hint that the grid for
the FFT does not have to be a rectangular one without elaborating this point.
In both cases it might be beneficial to consider a more general sampling,
i.e.~sampling on  anisotropic lattices.

Besides the theory of a discrete Fourier transform (DFT) on abelian groups,
also known as generalized Fourier transform~\cite{AhlanderMunthe-Kaas2005}, the
DFT has also been generalized to arbitrary sampling lattices, e.g.~in order
to derive periodic
wavelets~\cite{LangemannPrestin2010WaveletAnalysis,BergmannPrestin2014dlVP} and
a corresponding fast Fourier and fast wavelet transform~\cite{Bergmann2013FFT}.
The computational complexity on these lattices even stays the same as on the
usual rectangular or pixel grid.

Furthermore using the theory of rank-1-lattices, Kämmerer et
al.~\cite{KaemmererPottsVolkmer2015a,KaemmererPottsVolkmer2015b} and Potts and
Volkmer~\cite{PottsVolkmer2015} derive several adaptive schemes to approximate
both a certain set of frequencies as well as a set of sampling points and
derive approximation errors for functions of certain smoothness. This also
includes a constructive derivation of the vector that generates the lattice.
For these special lattices, the Fourier transform even in high-dimensional space
reduces to a one-dimensional Fourier transform and hence reducing both the
organization of the sampled data and the computational cost to compute the FFT.
The theory of rank-1-lattices therefore both allows for directly taking known
anisotropic properties of a function into account and thereby reducing the
necessary number of sample values or measurements by adapting the lattice.
Furthermore it also reduces the computational cost or data
organization overhead due to the reduction from a high-dimensional FFT to a
one-dimensional one.

In this paper we generalize the Basic Scheme by Moulinec and Suquet to arbitrary
anisotropic periodic lattices.
This introduces the possibility to prefer directions other than the coordinate
axes in the reference volume and hence in  the solution. This allows for
aligning the basis functions with the dominant orientations of the geometry and
controlling the refinement in these directions.
This generalization of the Basic Scheme to arbitrary anisotropic sampling
lattices introduces the form of the grid as a algorithmic parameter
without additional computational costs. For a special set of rank-1-lattices,
after sampling in a high-dimensional space, the computation of the fast Fourier
transform even reduces to a one-dimensional FFT.
Therefore, additionally to the new possibility of choosing directions of
preference, one can also choose these to reduce the computational efforts.

The remainder of the paper is organized as follows.
In Section~\ref{sec:Preliminaries} we establish the preliminaries regarding
the parametrization, properties of anisotropic lattices and their patterns on
the unit cube. Further, we introduce
the FFT on such patterns, where the usual tensor product grid is a special case.
Exemplary for a homogenization problem we introduce the periodic equations of
quasi-static elasticity in Section~\ref{sec:homogenization} and elaborate on the
unmodified Basic Scheme how to generalize it. Based on this we explain the
difference between making a coordinate transformation and choosing a lattice
adapted to the geometry of the lattice.
In Section~\ref{sec:geometries} we generalize two known geometries to an anisotropic setting: the laminate structure and the Hashin structure that serves as the main analytical example for this work.
Both are anisotropic structures within isotropic material laws that provide an
analytic solution for the strain field and the effective matrix.
This allows us to study effect of the pattern orientation on
the solution and the effective properties in Section~\ref{sec:numerics}.

%
%
\section{Preliminaries}\label{sec:Preliminaries}
Throughout this paper we will employ the following notation: The symbols \(a\in\mathbb C\),
\(\vect{a}\in\mathbb C^d\) and \(\mat{A}\in\mathbb C^{d\times d}\) denote scalars,
vectors, and matrices, respectively. The only exception from this
are~\(f,g,h\) which are reserved for functions.
We denote the inner product of two
vectors by \(\vect{a}^\tT\vect{b} \coloneqq\sum_i a_ib_i\) and reserve the
symbol \(\langle\cdot,\cdot\rangle\) for inner products of two functions or
two generalized sequences, respectively. For a complex number \(a = b+\im c\), \(b,c\in\mathbb R\), we denote the complex conjugate by \(\overline{a} \coloneqq b+\im c\).

Usually, we are concerned with \(d\)-dimensional data, where \(d=2,3\),
but the theory can also be written in arbitrary dimensions. 
Sets are denoted by capital case calligraphic letters, e.g.~\(\mathcal P(\cdot)\)
or~\(\mathcal G(\cdot)\) and the same for the Fourier transform~\(\mathcal F(\cdot)\);
all of these might depend on a scalar~\(n\) or matrix~\(\mat{M}\) given in brackets.
We denote second-order tensors by small Greek letters
as~\(\lambda,\epsilon\) with entries \(\lambda_{ij}\)
are indexed again by scalars \(i,j\) and similarly we denote fourth-order tensors
by capital calligraphic letters, where \(\mathcal C\) is the most prominent one.
Finally, constants like Euler's number~\(\e\) or the imaginary unit \(\im\), i.e.\ \(\im^2 = -1\),
are set upright.

\subsection{Arbitrary patterns and the Fourier transform}
\label{subsec:preliminariesFourier}
The space of functions we are concerned with is the Hilbert
space~$L^2(\T^d)$ of (equivalence classes of) square integrable functions
on the \(d\)-dimensional torus~$\T \cong [-\pi,\pi)^d$ with inner product
\begin{equation*}
	\langle f,g \rangle 
	=\frac{1}{(2\pi)^d}
		\int_{\T^d}
			f(\vect{x})\overline{g(\vect{x})}
		\dx\vect{x},
		\qquad
		f,g \in L^2(\T^d)
	\text{.}
\end{equation*}
In several cases, the functions of interest are tensor-valued.
For these functions, we take the tensor product of the Hilbert space,
e.g.~\(L^2(\mathbb T^d)^{n\times n}\) for the space of
functions~\(f\colon \mathbb T^d \to \mathbb C^{n\times n}\) that have values
being~\(n\times n\)-dimensional matrices. The following
Fourier transform can be generalized to these tensor product spaces by performing
the operations element wise. We restrict the following preliminaries of this
subsection therefore to the case of \(L^2(\mathbb T^d)\).

Every function $f \in L^2(\T^d)$ can be written in its Fourier series
representation
\begin{equation}\label{eq:fourier-series}
	f(\vect{x})
	= \sum_{\vect{k} \in \mathbb Z^d} c_{\vect{k}}(f)\e^{\im\vect{k}^\tT\vect{x}},
\end{equation}
introducing the the multivariate Fourier
coefficients~\(c_{\vect{k}}(f) = \langle f,\e^{\im\vect{k}^\tT\circ}\rangle\),
\(\vect{k}\in\mathbb Z^d\). The equality in \eqref{eq:fourier-series} is meant in  \(L^2(\T^d)\) sense.
We denote by $\vect{c}(f)
	= \bigl\{(c_{\vect{k}}(f)\bigr\}_{\vect{k}\in\mathbb Z^d}
	\in \ell^2(\mathbb Z^d)$ the generalized sequences which form a Hilbert
space with the inner product
\begin{equation*}
	\langle\vect{c},\vect{d}\rangle
		= \sum_{\vect{k} \in \mathbb Z^d}c_{\vect{k}}\overline{d_{\vect{k}}},
		\qquad \vect{c},\vect{d}\in \ell^2(\mathbb Z^d)\text{.}
\end{equation*}
The Parseval equation reads
\begin{equation}
		\langle f, g \rangle
		= \langle \vect{c}(f),\vect{c}(g) \rangle
		= \sum_{\vect{k} \in \mathbb Z^d} c_{\vect{k}}(f) \overline{c_{\vect{k}}(g)}\text{.}
		\label{eq:parseval}
\end{equation}

\paragraph{The pattern and the generating set}
For any regular matrix $\mat{M} \in \mathbb Z^{d\times d}$, we define the
congruence relation for $\vect{h},\vect{k} \in \mathbb Z^d$ with respect
to~$\mat{M}$ by
\begin{equation*}
		\vect{h} \equiv \vect{k} \bmod \mat{M}
		\Leftrightarrow \exists\,\vect{z} \in \mathbb Z^d\colon \vect{k} = \vect{h} + \mat{M}\vect{z}\text{.}
\end{equation*}
We define the lattice
\[
	\Lambda(\mat{M}) \coloneqq \mat{M}^{-1}\mathbb Z^d
	= \{\vect{y}\in\mathbb R^d : \mat{M}\vect{y} \in \mathbb Z^d\},
\]
and the pattern \(\Pattern(\mat{M})\) as any set of congruence representant of
the lattice with respect to \(\bmod\ 1\), e.g.\ $ \Lambda(\mat{M})\cap[0,1)^d$
or~$\Lambda(\mat{M})\cap\bigl[-\tfrac{1}{2},\tfrac{1}{2}\bigr)^d$. For the rest
of the paper we will refer to the set of congruence class representants in the 
symmetric unit cube \(\bigl[-\tfrac{1}{2},\tfrac{1}{2}\bigr)^d\). The generating
set \(\generatingSet(\mat{M})\) is defined by \(\generatingSet(\mat{M}) 
\coloneqq \mat{M}\Pattern(\mat{M})\) for any pattern \(\Pattern(\mat{M})\). For
both, the number of elements is given by \(
	\abs{\Pattern(\mat{M})}
	=\abs{\generatingSet(\mat{M})}
	=\abs{\det{\mat{M}}}
	\eqqcolon m,
\) which follows directly from~\cite[Lemma II.7]{deBoorHoelligRiemenscheider1993BoxSplines}.

Finally for any factorization \(\mat{M} = \mat{J}\mat{N}\) of an integer matrix
\(\mat{M}\in\mathbb Z^{d\times d}\) into two integer matrices \(\mat{J},\mat{N}\in\mathbb Z^d\),
we have
\[
	\vect{x}\in\Lambda(\mat{N}) \Rightarrow \mat{N}\vect{x} \in\mathbb Z^d
	\Rightarrow \mat{J}\mat{N}\vect{x}\in\mathbb Z^d,
\]
and hence \(\Lambda(\mat{N})\subset\Lambda(\mat{M})\). By construction of
the pattern, we directly obtain \(\Pattern(\mat{N}) \subset \Pattern(\mat{M})\);
see~\cite[Lemma 2.7]{LangemannPrestin2010WaveletAnalysis}
and~\cite[Section 2]{Bergmann2013FFT} for a more general introduction.
We call the smaller pattern~\(\Pattern(\mat{N})\) a \emph{subpattern}
of~\(\Pattern(\mat{M})\).
Note that this is not commutative with respect to \(\mat{J}\) and \(\mat{N}\).
Looking at the generating sets for the decomposition \(\mat{M} = \mat{JN}\) we have
\(\generatingSet(\mat{N}^\tT) \subset \generatingSet(\mat{M}^\tT)\). A subpattern of a tensor product grid, the so called quincunx pattern,
is shown in Fig.~\ref{fig:patterns}, left.
%
%
\subsection{A fast Fourier transform on patterns}\label{subsec:FFT}
The discrete Fourier transform on the pattern $\Pattern(\mat{M})$ is defined~\cite{ChuiLi:1994} by
\begin{equation}\label{eq:Fouriermatrix}
	\mathcal F(\mat{M})
	\coloneqq
	\frac{1}{m}
	\Bigl(
		\e^{- 2\pi \im \vect{h}^\tT\vect{y}}
	\Bigr)_{%
		\vect{h} \in \generatingSet(\mat{M}^\tT),\,%
		\vect{y} \in \Pattern(\mat{M})%
},
\end{equation}
where $\vect{h}\in\generatingSet(\mat{M}^\tT)$ indicate the rows and
$\vect{y} \in \Pattern(\mat{M})$ indicate the columns of the Fourier
matrix~\(\mathcal F(\mat{M})\). This discrete Fourier transform was also
investigated in~\cite{Bergmann2013FFT,LangemannPrestin2010WaveletAnalysis}.
For both the pattern \(\Pattern(\mat{M})\) and the generating
set~\(\generatingSet(\mat{M}^\tT)\), an arbitrary but fixed ordering has to be chosen.
The discrete Fourier transform 
on $\Pattern(\mat{M})$ is defined for a
vector~$\vect{a} = (a_{\vect{y}})_{\vect{y}\in\Pattern(\mat{M})}\in\mathbb C^m$
arranged in the same ordering as the columns in~\eqref{eq:Fouriermatrix} by
\begin{equation}\label{eq:FourierTransform}
	\vect{\hat a} = (\hat a_{\vect{h}})_{\vect{h}\in\generatingSet(\mat{M}^\tT)}
	= \mathcal F(\mat{M})\vect{a},
\end{equation}
where the resulting vector $\vect{\hat a}$ is ordered as the columns
of~$\mathcal F(\mat{M})$ in~\eqref{eq:Fouriermatrix}. 

For a diagonal
matrix~\(\mat{M} = \operatorname{diag}(n,n)\in\mathbb N^{2\times 2}\) having the
same entry \(n\in\mathbb N\) on both diagonal entries,
the pattern~\(\Pattern(\mat{M})\) is the
set \(\bigl(y_1,y_2)^\tT, y_1,y_2\in\frac{1}{m}\{0,\ldots,m-1\}\).
The generating set~\(\generatingSet(\mat{M})\) then
reads~\(\{ \vect{k}\in\mathbb Z^2 : 0\leq k_j \leq m-1, j=1,2\}\).
Both have \(m = n^2\) elements.
The Fourier transform for this diagonal matrix is just the usual 2D DFT.

\begin{figure}[tbp]
	\begin{subfigure}[t]{.32\textwidth}\centering
		\includegraphics{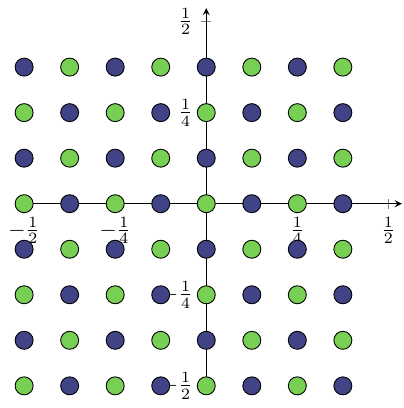}
	\end{subfigure}
	\begin{subfigure}[t]{.32\textwidth}\centering
		\includegraphics{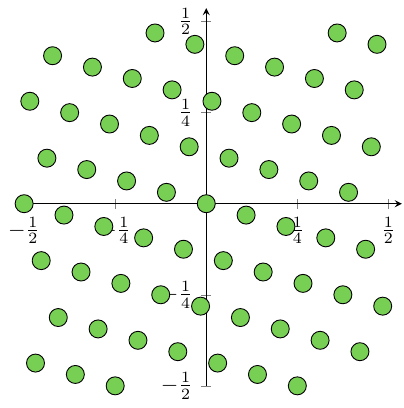}
	\end{subfigure}
	\begin{subfigure}[t]{.32\textwidth}\centering
		\includegraphics{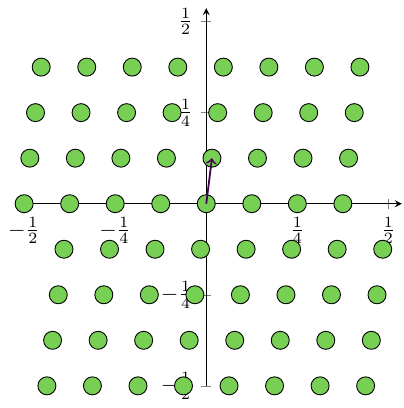}
	\end{subfigure}
	\caption[Three Patterns.]{
	By choosing \(\mat{M}_1 = \mat{JN} = \bigl( \begin{smallmatrix}1&-1\\1&1\end{smallmatrix}\bigr)\bigl(\begin{smallmatrix}4&4\\-4&4\end{smallmatrix}\bigr)\) we obtain a usual rectangular grid pattern \(\Pattern(\mat{M}_1)\) (left) of a diagonal matrix and its sub pattern \(\Pattern(\mat{N}\) (dark), a quincux pattern. We can prefer certain directions like for pattern~\(\Pattern(\mat{M}_2)\), \(\mat{M}_2= \bigl(\begin{smallmatrix} 8&-4\\2&7\end{smallmatrix}\bigr)\), (middle). Certain patterns, like \(\Pattern(\mat{M}_3)\), \(\mat{M}_3= \bigl(\begin{smallmatrix}8&-1\\0&8\end{smallmatrix}\bigr)\), (right) are even generated by only one generating vector. Note that all matrices have the same determinant and hence the patterns have the same number of points.}
	\label{fig:patterns}
\end{figure}
Fig.~\ref{fig:patterns} illustrates that the Fourier
transform \eqref{eq:FourierTransform} generalizes the usual discretization on a
pixel grid and enables to prefer certain directions by choosing different patterns having the same number of points.

To efficiently implement the Fourier
transform~\eqref{eq:FourierTransform} on an arbitrary
pattern~\(\Pattern(\mat{M)}\) for a regular integer
matrix~\(\mat{M}\in\mathbb Z^{d\times d}\), we have to fix a certain order of
the elements therein.
Following the construction in~\cite{Bergmann2013FFT}, we use the
Smith normal form \(\mat{M} = \mat{QER}\), where \(\mat{Q},\mat{R}\) are of
determinant \(1\) and~\(\mat{E} = \diag(e_1,\ldots,e_d)\)
is the diagonal matrix of elementary divisors where \(e_j\) is a
divisor of~\(e_{j+1}\), \(j=1,\ldots,d-1\).
We further denote by~\(d_{\mat{M}}
	\coloneqq \lvert \{j : e_j>1\}\rvert\)
the \emph{dimension of the pattern}.
For the special case, that \(d_{\mat{M}}=1\), the lattice is also called
\emph{rank-\(1\)-lattice}. Such a lattice is shown in Fig.~\ref{fig:patterns}
~(right)

Introducing the \emph{pattern basis vector(s)}
\begin{equation}\label{eq:patternBasisV}
	\vect{y}_j \coloneqq
	\frac{1}{e_{d+d_{\mat{M}}+j}}\vect{e}_{d-d_{\mat{M}}+j},\qquad j=1,\ldots,d_{\mat{M}},
\end{equation}
where \(\vect{e}_j\) denotes the \(j\)th unit vector, we obtain a basis for the
pattern. Hence we can write
\begin{equation*}
\vect{y} = \sum_{j=1}^{d_{\mat{M}}}\lambda_j\vect{y}_j,\qquad
\lambda_j\in\{0,\ldots,e_j-1\},\ j=1,\ldots,d_{\mat{M}},
\end{equation*}
where the summation is meant on the congruence classes,
i.e.~with respect to \(\bmod\ 1\) onto the pattern~\(\Pattern(\mat{M})\), we obtain a
unique addressing for each pattern point using the coefficients \(\lambda_1,\ldots,\lambda_{d_{\mat{M}}}\).
With the lexicographical ordering of the vectors \((0,\ldots,0)^\tT,\ldots,\)\\
\((\lambda_1-1,\ldots,\lambda_{d_{\mat{M}}}-1)\)
one not only obtains an array representation of any coefficient vector
\[
	\vect{a} = \bigl(a_{\vect{y}}\bigr)_{\vect{y}\in\Pattern(\mat{M})}
	= \bigl(a_{\lambda_1,\ldots,\lambda_{d_{\mat{M}}}}\bigr)_{\lambda_1=0,\ldots,\lambda_{d_{\mat{M}}}=0}^{e_{d-d_{\mat{M}}+1}-1,\ldots,e_d-1},
\]
but similarly also for any vector \(\hat{\vect{a}}\) corresponding to the
generating set using the \emph{generating set basis vector(s)}
\begin{equation}\label{eq:genSetBasisV}
	\vect{h}_j \coloneqq\mat{M}^\tT\tilde{\vect{y}}_j = \mat{R}^\tT\vect{e}_j,\qquad j=1,\ldots,d_{\mat{M}},
\end{equation}
where \(\tilde{\vect{y}}\) denotes the basis vector(s) of \(\Pattern(\mat{M}^\tT)\) constructed as above.
Note that the pattern dimension \(d_{\mat{M}}=d_{\mat{M}^\tT}\) and the elementary
divisors~\(e_j\) are identical for the patterns~\(\Pattern(\mat{M})\) and \(\Pattern(\mat{M}^\tT)\).

With these fixed orderings of the vector entries of \(\vect{a}\) and \(\hat{\vect{a}}\),
the Fourier transform~\eqref{eq:FourierTransform} can be computed using
an ordinary \(d_{\mat{M}}\)-dimensional Fourier transform even having the
same leading coefficient in its complexity of \(\mathcal O(m\log m)\)~\cite[Theorem~2]{Bergmann2013FFT}. Note that for rank-1-lattices, the Fourier transform on the pattern even reduces
to a one-dimensional FFT for patterns in 2, 3 or even more dimensions.
%
\paragraph{Sampling and aliasing}
Let \(f\in L_2(\mathbb T^d)\) denote a square integrable function on the torus,
such that its Fourier series \eqref{eq:fourier-series} converges absolutely, i.e.
\[
	\sum_{\vect{k}\in\mathbb Z^d} \abs{c_{\vect{k}}((f))} < \infty.
\]
Sampling \(f\) at the points given by pattern \(\Pattern(\mat{M})\) of a regular matrix \(\mat{M}\in\mathbb Z^{d\times d}\), i.e.~\(a_{\vect{y}} \coloneqq f(2\pi\vect{y})\), \(\vect{y}\in\Pattern(\mat{M})\),
and performing a discrete Fourier transform
we obtain the \emph{discrete Fourier coefficients} \(c_{\vect{h}}^{\mat{M}}(f) = \hat a_{\vect{h}}\), \(\vect{h}\in\generatingSet(\mat{M}^\tT)\), where
\(\vect{\hat a} = \mathcal F(\mat{M})\vect{a}\).
A relation between the Fourier coefficients \(c_{\vect{k}}(f)\) and the
discrete Fourier coefficients~\(c_{\vect{h}}^{\mat{M}}(f)\)
is given by the following Lemma.
\begin{lemma}\label{lem:Aliasing}
	Let \( f\in L^2(\T^d) \) with absolutely convergent Fourier series and
	the regular matrix~\(\mat{M}\in\mathbb Z^{d\times d}\) be given.
	Then the discrete Fourier coefficients~\( c_{\vect{h}}^{\mat{M}}(f) \),
	\(\vect{h}\in\generatingSet(\mat{M}^\tT)\), fulfill the relation
	\begin{equation}\label{eq:aliasing}
		c_{\vect{h}}^{\mat{M}}(f) 
		= 
		\sum_{\vect{z}\in\mathbb Z^d} c_{\vect{h}+\mat{M}^\tT\vect{z}}(f)
		,\qquad \vect{h}\in\generatingSet(\mat{M)^\tT}\text{.}
	\end{equation}
\end{lemma}
For a proof, see~\cite[Lemma~2]{BergmannPrestin2014Interpolation}. This Lemma
is also called \emph{Aliasing formula} and can be interpreted as follows: If the
Fourier coefficients~\(c_{\vect{k}}(f)\) of \(f\) decay slowly along a certain
direction~\(\vect{h}_j\) being a basis vector of \(\generatingSet(\mat{M}^\tT\))
and the corresponding \(\epsilon_j\) is small, then the effect of the summands
\(\vect{z}\neq \vect{0}\) is quite large or in other words the approximation
\(c_\vect{h}^{\mat{M}}(f) \approx c_{\vect{h}}(f)\) is not sufficient enough.
This might be e.g.~due to presence of an edge orthogonal to \(\vect{h}_j\). If on the other
hand, the spectrum is bounded in this direction and \(\epsilon_j\) is large enough,
then the approximation is of better quality.
%
\paragraph{Pattern congruence classes}
Following the pattern classification, cf.~\cite[Section 2.4]{LangemannPrestin2010WaveletAnalysis} we notice that \(\Pattern(\mat{M}) = \Pattern(\mat{N})\) whenever \(\mat{M} = \mat{Q}\mat{N}\) holds for a matrix \(\mat{Q}\in\mathbb Z^{d\times d}\), \(\lvert\det\mat{Q}\rvert = 1\). We define the \emph{pattern congruence} of two matrices by \(\mat{M}\sim_{\Pattern}\mat{N}\) whenever they generate the same pattern. By~\cite[Lemma 2.5]{LangemannPrestin2010WaveletAnalysis} there exists a congruence class representant \(\mat{M}^{\circ} = (m^{\circ}_{i,j})_{j,i=1}^{d,d}\) in every congruence class, such that \(\mat{M}^{\circ}\) is upper triangular, and \(0\leq m^{\circ}_{i,j} < m^{\circ}_{j,j}\) for all \(i<j\).

An easy consequence of this is, that for a diagonal matrix \(\mat{M} = \bigl(
	\begin{smallmatrix}
	m_1 & 0 \\0 & m_2
\end{smallmatrix}\bigr)\) choosing \(\mat{Q} = \bigl( \begin{smallmatrix}
	1&n\\0&1
\end{smallmatrix})\), \(n\in\mathbb N\), reveals, that all matrices \(\mat{M}_n \coloneqq \mat{Q}\mat{M} = \bigl(
	\begin{smallmatrix}
		m_1 & m_2n\\0&m_2
	\end{smallmatrix}
\bigr)\) possess the same pattern or in other words the same sampling points on the torus. However, their corresponding generating sets \(\generatingSet(\mat{M}_n^\tT)\) differ. This can be interpreted as choosing different directional sine and cosine functions that can be defined on the same set of points in order to analyze given discrete data or employing the periodicity of \(c_{\vect{h}}^{\mat{M}}(f)\) with respect to \(\mat{M}^\tT\) that is implied by~\eqref{eq:aliasing}. This introduces the possibility of anisotropic analysis and interpretation even on a usual pixel grid.

From the computational point, two different matrices
\(\mat{M}_1 \sim_{\Pattern} \mat{M}_2\)
having the same pattern \(\Pattern(\mat{M}_1) = \Pattern(\mat{M}_2)\)
only require an rearranging of the data arrays due to the change of the basis
vectors, i.e.~using \(\vect{y}^{(1)}_j\) or \(\vect{y}^{(2)}_j\),
\(j=1,\ldots,d_{\mat{M}_1}=d_{\mat{M}_2}\), when addressing \(\vect{a}\)
in~\eqref{eq:FourierTransform} and similarly the ordering with respect to the
generating sets~\(\generatingSet(\mat{M}^\tT_i)\), \(i=1,2\). This rearrangement can easily be computed, see~\cite[Section 5.2]{Bergmann2013FFT}.

\begin{figure}\centering
	\begin{subfigure}[t]{.32\textwidth}\centering
		\includegraphics{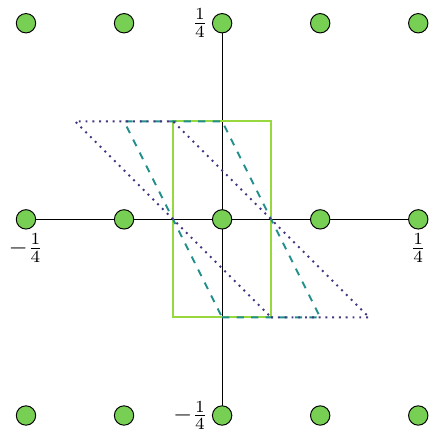}
	\end{subfigure}
	\hspace{.16\textwidth}
	\begin{subfigure}[t]{.32\textwidth}\centering
		\includegraphics{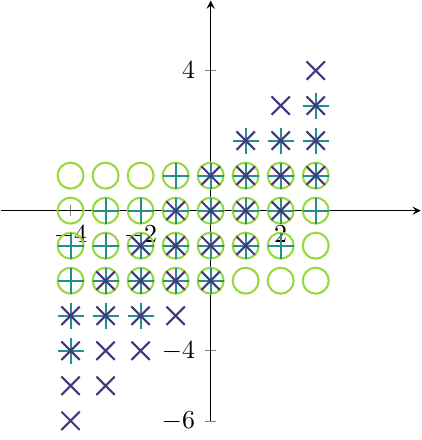}
	\end{subfigure}
	\caption[]{ The three matrices~\(\mat{M} = \bigl(
		\begin{smallmatrix}
			8 & 4n\\0&4
		\end{smallmatrix}
	\bigr), n=0,1,2\) induce the same sample point set but different scaled unit cells \(\frac{1}{2}\mat{M}_n^{-1}[-1,1]^2\) (solid, dashed, dotted; left) and different generating sets \(\generatingSet(\mat{M}_n^\tT)\) (\(\circ\), \(+\), \(\times\); right, respectively) yielding different frequency sets for the Fourier transform.}\label{fig:DiffSample}
\end{figure}
\begin{example}
	For the simple matrix \(\mat{M} = \mat{M}_0 = \bigl(\begin{smallmatrix}
		8&0\\0&4
	\end{smallmatrix}\bigr)\) the pattern is just the rectangular sampling grid,
	a subset of which is displayed in Fig.~\ref{fig:DiffSample}
	~(left). The two matrices \(\mat{M}_1 = \bigl(
	\begin{smallmatrix}
		8&4\\0&4
	\end{smallmatrix}
	\bigr)\) and \(\mat{M}_2 = \bigl(
	\begin{smallmatrix}
		8&8\\0&4
	\end{smallmatrix}
	\bigr)\) possess the same sample points. However, their scaled unit cells \(\frac{1}{2}\mat{M}_i^{-1}[-1,1]^2\), also shown in Fig.~\ref{fig:DiffSample}
	~(left), differ. This illustrates the different directional preference of these matrices, though they possess the same sampling lattice.  While the diagonal matrix resembles the form of (stretched) pixels, the shear introduced by looking at other sine/cosine terms can be clearly seen for \(i=1,2\) for the dashed and dotted parallelotopes.
	The three different generating sets \(\generatingSet(\mat{M}_n^\tT)\), \(n=1,2,3\) are denoted as circles, pluses and crosses in Fig.~\ref{fig:DiffSample}
	(right), respectively.
	Hence both figures illustrate one way to visualize the anisotropy.
\end{example}
%
%
%
\section{Homogenization}\label{sec:homogenization}
We want to investigate composite structures which consist of a finite number of
materials composed into one material. Common examples are fiber reinforced
polymers~\cite{KMS:2015Homogenization}, polycrystalline
structures~\cite{Lebensohn} or metall foams~\cite{Liebscher2014}. As these
structures are typically small they can seldom be resolved exactly in
simulations of larger structures. This gives rise to the concept of
homogenization where the composite material is replaced by a homogeneous one
having the same relevant, i.e.~macroscopic, properties.
The basic assumption to do this is that the microstructure is periodic, i.e.~it is
sufficient to look at a representative volume element (RVE) and that the scales
of the microstructure and macroscopic structure are separated.

These assumptions allow to calculate the homogenized behavior of the material,
the so called effective properties or effective matrix. This can be
inserted into a macroscopic calculation or can be used to determine the isotropy
of the structure or relevant elastic properties. Typical examples
of such problems involve the steady-state heat equation or the quasi-static
equation of linear elasticity which we want to focus on henceforth.

\subsection{The equation of quasi-static linear elasticity in homogenization}
\label{pg:Elasticity}

The partial differential equation (PDE) we consider as an example in this
publication is the equation of quasi-static linear elasticity. It
will serve as the basis to develop a numerical algorithm solving the PDE
later on.

Consider a periodic stiffness
distribution~$\Stiffness \in L^\infty(\T^d)^{d \times d \times d \times d}$
that is essentially bounded with major and minor symmetries, i.e.~%
with~$
\Stiffness_{ijkl}
	= \Stiffness_{jikl} = \Stiffness_{ijlk} =  \Stiffness_{klij}
$ characterizing the microstructure in the representative
volume element.
The entries in~$\Stiffness$ specify the material behavior,
e.g.\ for an isotropic material we have~$
	\Stiffness_{ijkl}
	= \lambda \delta_{ij}\delta_{kl} + \mu(\delta_{ik}\delta_{jl}
		+ \delta_{il}\delta_{jk})$
where $\lambda \in \R$ and $\mu \geq 0$ are the first and second Lam\'e
parameter, respectively.

For the variational formulation of the steady-state heat equation Vord{\v{r}}ejc
et.al.~\cite{Vondrejc2014} derive an equivalent integral formula. This approach
can be directly applied to the quasi-static linear elasticity equation as
follows. We define the space 
\begin{equation*}
	\mathcal{E}(\T^d) \coloneqq
\Bigl\lbrace \vect{v} \in L^2(\T^d)^{d \times d} : \vect{v} = \GradSym
\vect{u}, \vect{u} \in H^1(\T^d)^d \Bigr\rbrace
\end{equation*}
where
$\GradSym \vect{u} \coloneqq \frac12 (\Grad \vect{u} + (\Grad \vect{u})^\tT)$ is the symmetric
gradient operator applied to displacement $\vect{u}$ and $H^1(\T^d)$ is the Sobolev 
space with Sobolev index 1, i.e.~the space of functions from $L^2(\T^d)$
where the first weak derivative is also in $L^2(\T^d)$.

In the following we will use Einstein's summation convention. For the product between a symmetric tensor $\mathcal A \in \R^{d \times d \times d \times d}$ of fourth order
and a symmetric second-order tensor $\alpha \in \R^{d \times d}$ we introduce the notation
\begin{equation*}
\mathcal A \tensorProd \alpha \coloneqq (\mathcal{A}_{ijkl} \alpha_{kl})_{ij}.
\end{equation*}

\begin{definition}
The partial differential equation of quasi-static linear elasticity in
homogenization reads:

Find for a macroscopic strain $\epsilon^0 \in \R^{d \times d}$ the strain
function $\tilde{\epsilon} = \tilde{\epsilon}_{\epsilon^0} \in \mathcal{E}(\T^d)$ such that for all
$\nu \in \mathcal{E}(\T^d)$ 
\begin{equation}\label{eq:linear-elasticity-weak}
	\langle \nu, \Stiffness \tensorProd (\epsilon^0 + \tilde{\epsilon})
	\rangle = 0
\end{equation}
holds true.
\end{definition}
We further
call~$\Stiffness^\eff \in \R^{d \times d \times d \times d}$ the effective
matrix which is connected with the PDE above by
\begin{equation}
\Stiffness^\eff \tensorProd \epsilon^0 \coloneqq \int_{\T^d} \Stiffness
\tensorProd \bigl(\tilde{\epsilon} + \epsilon^0 \bigr) \dx \vect{x}.
\end{equation}

By \cite{Vondrejc2014} this is equivalent to solving an integral equation for
the strain $\epsilon = \tilde{\epsilon} + \epsilon^0$ introducing a constant
non-zero reference stiffness $\Stiffness^0$.

\begin{definition}
The Lippmann-Schwinger equation is given as:

Find $\epsilon \in L^2(\T^d)^{d \times d}$ such that
\begin{equation}
	\epsilon = \epsilon^0 - \GradSym (\Div \tensorProd \Stiffness^0 \tensorProd \GradSym)^{-1} \Div (\Stiffness - \Stiffness^0) \tensorProd \epsilon
\label{eq:lippmann-schwinger-weak}
\end{equation}
in weak sense.
\end{definition}

The divergence operator is formally the negative adjoint of the symmetric
gradient operator, i.e. $\Div = -\GradSym^\hH$, and
therefore~\eqref{eq:lippmann-schwinger-weak} reduces to finding a strain $\epsilon \in
L^2(\T^d)^{d \times d}$ such that
\begin{equation}
	\epsilon = \epsilon^0 - \GradSym (\GradSym^\hH \tensorProd \Stiffness^0
	\tensorProd \GradSym)^{-1} \GradSym^\hH (\Stiffness - \Stiffness^0)
	\tensorProd \epsilon,
\label{eq:lippmann-schwinger}
\end{equation}
see \cite{Schneider2016}.
The strain $\epsilon \in L^2(\T^d)^{d \times d}$ can be represented by its
Fourier series and we obtain for \(\vect{k} \in \Z^d\setminus\{\vect{0}\}\) that
\begin{equation}
	c_\vect{k}(\epsilon) =
	- \GradSym_\vect{k} (\GradSym^\hH_\vect{k} \tensorProd \Stiffness^0 
	\tensorProd \GradSym_\vect{k})^{-1} \GradSym^\hH_\vect{k}
		c_\vect{k} \bigl(
			 (\Stiffness - \Stiffness^0) \tensorProd \epsilon
		\bigr),
\label{eq:lippmann-schwinger-fourier}
\end{equation}
and \(c_{\vect{0}}(\epsilon) = \epsilon^0\).
The Fourier multiplier $\GradSym_\vect{k}$ hereby represents the action of the
operator~$\GradSym$ with respect to a Fourier coefficient index
\(\vect{k}\in\mathbb Z^d\backslash\{\vect{0}\}\), respectively. For $u \in
H^1(\T^d)^d$ it can be derived as
\begin{equation}
	c_{\vect{k}}(\GradSym \vect{u})
	\coloneqq
	\GradSym_\vect{k} c_\vect{k}(\vect{u})
	=	\frac{\im}{2} (\vect{k} c_\vect{k}(\vect{u})^\tT 
	+ c_\vect{k}(\vect{u}) \vect{k}^\tT).
\end{equation}

\subsection{Homogenization on anisotropic lattices}
The approach of Moulinec and Suquet~\cite{MoulinecSuquet1994,MoulinecSuquet1998}
to discretize~\eqref{eq:lippmann-schwinger-fourier} is based on collocation on a
Cartesian grid, see also~\cite{ZVNM:2010NumHom}. To take into account preferred directions in
composite we generalize this approach to arbitrary patterns
\(\Pattern(\mat{M})\) introduced in the Section~\ref{sec:Preliminaries}.

Let a regular integer matrix \(\mat{M}\in \Z^{d \times d}\) be given. Following
the idea of Moulinec and Suquet we collocate~\eqref{eq:lippmann-schwinger} at
the points~\(2\pi\vect{y}\in\mathbb T^d, \vect{y}\in\Pattern(\mat{M})\), of the
pattern. This discretization leads to
the problem of finding symmetric matrices $\epsilon_\vect{y} \in \R^{d \times
d}$ for each $\vect{y} \in \Pattern(\mat{M})$ such that for
\(\vect{y}\in\Pattern(\mat{M})\setminus\{\vect{0}\}\) we have
\begin{equation}
	\begin{split}
\epsilon_\vect{y} = -&\sum_{\vect{h} \in \generatingSet(\mat{M}^\tT)} 
\GradSym_\vect{h} \bigl( \GradSym^\hH_\vect{h} \tensorProd \Stiffness^0
\tensorProd \GradSym_\vect{h} \bigr)^{-1} \GradSym_\vect{h}\\
&\quad \times \frac{1}{m}\sum_{\vect{z} \in \Pattern(\mat{M})} (\Stiffness_\vect{z} - \Stiffness^0) \tensorProd
\epsilon_\vect{z} \e^{-2 \pi \im \vect{h}^\tT \vect{z}} \e^{2 \pi \im \vect{h}^\tT \vect{y}}
	\end{split}
\end{equation}
and $\epsilon_{\vect{0}} = \epsilon^0$.
This discretization gives rise to the generalized basic scheme for patterns based
on~\cite{MoulinecSuquet1998} summarized in Algorithm~\ref{alg:fixed-point}.
\begin{algorithm}[t]
\begin{algorithmic}
\State $\epsilon^{(0)}_\vect{y} \gets \epsilon^0$ for all $\vect{y} \in \Pattern(\mat{M})$
\State $n \gets 0$
\Repeat
\State $\tau^{(n+1)}_\vect{y}
		\gets \left(\Stiffness_\vect{y}
			- \Stiffness^0\right) \tensorProd \epsilon_\vect{y}^{(n)},\quad \vect{y}\in\Pattern(\mat{M})$
\State $\hat{\tau}^{(n+1)} \gets \Fourier(\mat{M}) \tau^{(n+1)}$
\State $\hat{\epsilon}_\vect{h}^{(n+1)}\gets - \GradSym_\vect{h} 
\bigl( \GradSym^\hH_\vect{h} \tensorProd \Stiffness^0
\tensorProd \GradSym_\vect{h} \bigr)^{-1} \GradSym_\vect{h} 
\hat{\tau}_\vect{h}^{(n+1)},\quad\vect{h}\in\generatingSet(\mat{M}^\tT)\setminus\{\vect{0}\}$
	\State $\hat{\epsilon}_\vect{0}^{(n+1)} \gets \epsilon^0 $
\State $\epsilon^{(n+1)} \gets \Fourier^{-1}(\mat{M}) \hat\epsilon^{(n+1)}$
\State $n \gets n+1$
\Until{a convergence criterion is reached}
\end{algorithmic}
\caption{Fixed-point algorithm on patterns.}
\label{alg:fixed-point}
\end{algorithm}

\subsection{Anisotropic unit cells}

By~\cite[Theorem 1.8]{Bergmann2013Thesis}, the necessary structure for the
sampling pattern $\Pattern(\mat{M})$ is an additive group structure. All these groups can be characterized by the congruence class representants \(\mat{M}^{\circ}\). Despite from that approach, another one is as follows: let
$\mat{L} \in \R^{d \times d}$ be a regular matrix and $\mathcal N \subset \R^d$
a set such that $(\mathcal N, + \mod \mat{L})$ is a group endowed with the
addition $+ \mod \mat{L}$. Then one can find a corresponding pattern
$\Pattern(\mat{M})$ such that $\mat{L} \Pattern(\mat{M}) = \mathcal N$, i.e.~we
can define a discrete Fourier transform on such a group.

\begin{figure*}
	\begin{subfigure}[b]{.24\textwidth}\centering
		\includegraphics{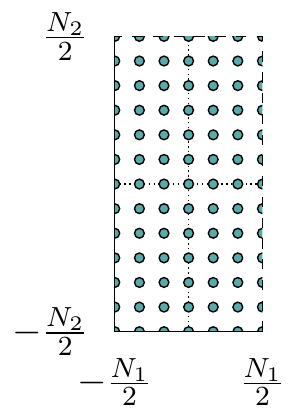}
	\end{subfigure}
	\begin{subfigure}[b]{.24\textwidth}\centering
		\includegraphics{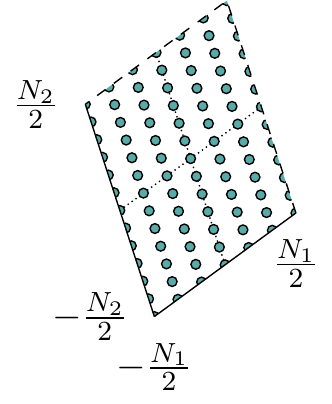}
	\end{subfigure}
	\begin{subfigure}[b]{.24\textwidth}\centering
		\includegraphics{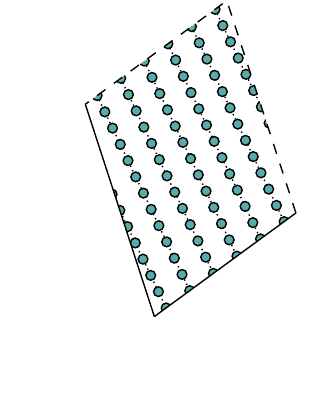}
	\end{subfigure}
	\begin{subfigure}[b]{.24\textwidth}\centering
		\includegraphics{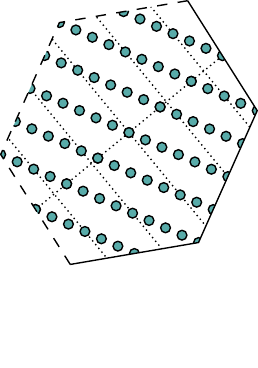}
	\end{subfigure}
	\caption[]{While the usual of a diagonal matrix indtroduces a rectangular grid (most left), the matrix \(\mat{L}\) of the transformed pattern \(\Pattern_{\mat{L}}(\mat{M})\) introduces rotation and scaling of the unit cube (left). Furthermore, in such a cell, a rank-1-lattice can be used (right), and it is even possible to use arbitrary shapes like an arbitrary hexagonal shape \(\mathcal S\) (most right).
	Dotted lines indicate lines of the basis vectors and all pattern consist of 72 points, whose boundary elements are repeated on the other side (the half circles each).
	}\label{fig:unitcells}
\end{figure*}

This setting allows for a more generalized notion of a pattern with examples
given in Figures~\ref{fig:unitcells}.

\begin{definition}\label{def:L-pattern}
Let $\mat{M} \in \Z^{d \times d}$ be an regular integral matrix and $\mat{L} \in \R^{d \times
d}$ be a regular matrix. Then the \emph{transformed pattern} is defined by
\begin{equation*}
\Pattern_\mat{L}(\mat{M}) \coloneqq \mat{L} \Lambda(\mat{M}) \cap \mat{L}
\bigl[-\tfrac12,\tfrac12\bigr)^d.
\end{equation*}
\end{definition}

One can even take \emph{any} set of integer points inside a certain shape
\(\mathcal S\), where the shifts~\(\mathcal S + \mat{L}\vect{z}\), tile the
\(\mathbb R^d\), i.e.~for all \(\vect{y}\in\mathbb R^d\) exists a
unique~\(\vect{z}\in \mathbb Z^d\) such that~\(\vect{y}-\vect{z}\in\mathcal S\),
e.g see Fig.~\ref{fig:unitcells} (most right)

This gives rise to a huge variety of unit
cells to model the microscopic periodic media.
The transformed patterns \(\Pattern_\mat{L}(\mat{M})\) especially introduce the
possibility to take anisotropic cell structures inside the unit cell or RVE into account.

Definition~\ref{def:L-pattern} can also be interpreted in terms of a coordinate
transformation. 
Consider a regular matrix $(A_{ij})_{i,j} = \mat{A} \in \R^{d
\times d}$ and transformed coordinates $\tilde{\vect{x}} = \mat{A} \vect{x}$
with $\vect{x} \in [-\frac12,\frac12)^d$. Let $\vect{u}(\vect{x})$ solve the PDE
\eqref{eq:linear-elasticity-weak} with stiffness distribution
$\Stiffness(\vect{x})$ and macroscopic strain $\epsilon^0$.

By \cite[Section 8.3]{Milton2002} the displacement
$\tilde{\vect{u}}(\tilde{\vect{x}})$ that solves
\eqref{eq:linear-elasticity-weak} with transformed
\begin{align*}
\tilde{\Stiffness}_{ijkl}(\tilde{\vect{x}}) &= A_{im} A_{jn}
A_{ko}A_{lp} \Stiffness_{mnop}(\vect{x}),
\end{align*}
and \(\tilde{\epsilon^0} = \mat{A}^{-T} \epsilon^0 \mat{A}^{-1}\)
is connected to $\vect{u}(\vect{x})$ by $\tilde{\vect{u}}(\tilde{\vect{x}}) =
\mat{A}^{-\tT} \vect{u}(\vect{x})$. Discretizing $\vect{u}(2 \pi \vect{y})$ with $\vect{y} \in
\Pattern(\mat{M})$ leads to
\begin{align*}
\vect{u}(2 \pi \vect{y}) = \mat{A}^\tT \tilde{\vect{u}}(2 \pi \mat{A} \vect{y})
= \mat{A}^\tT \tilde{\vect{u}}(2 \pi \tilde{\vect{y}})
\end{align*}
with \(\tilde{\vect{y}} = \mat{A} \vect{y} \in \Pattern(\mat{M})\) or equivalently \(\tilde{\vect{y}} \in \Pattern_\mat{A}(\mat{M})\).

Compared to the usual coordinate transform,
c.f.~\cite{Lahellec2003,Rietbergen1996}, to the unit cell, the homogenization on
patterns allows for further preference of directions in the cell other than the
transformed coordinate axes.

\subsection{Convergence of the discretization}
Let \(n\) be the smallest eigenvalue
of \(\mat{M}\) being larger than \(1\),
where~\(\mat{M}\in\mathbb Z^{d\times d}\). Then interpolation error on the
generating set~$\generatingSet(\mat{M}^\tT)$ can be bounded from above by the
diagonal matrix \(\tilde{\mat{M}} \coloneqq \operatorname{diag}(
	\lfloor n\rfloor,\ldots,\lfloor n\rfloor)\in\mathbb Z^{d\times d}\),
because the hypercube \(\frac{1}{2}\tilde{\mat{M}}[-1,1]^d\) is contained in the
parallelepiped \(\tfrac{1}{2}\mat{M}^\tT[-1,1]^d\) which is used to
define~$\generatingSet(\mat{M}^\tT)$ and hence $\generatingSet(\tilde{\mat{M}}^\tT)\subset\generatingSet(\mat{M}^\tT)$. From any sequence of discretizations~\(\mat{M}_i\) whose determinants \(m_i\) tend to infinity, we also obtain, that the smallest eigenvalue tends to infinity because the matrix has to span the frequency lattice \(\mathbb Z^d\). With this argument the convergence of the such discretization sequence is
dominated by a standard case of a Cartesian grid and the convergence proof
of e.g.~\cite{Schneider2015Convergence} the convergence in this setting follows directly.

\section{Geometries and analytic solutions}\label{sec:geometries}
To examine the effects the pattern matrix has on the solution,
this section describes two problems where the effective stiffness tensor and
an analytic expression for the strain field are known.
The first structure is a laminate which as a basically one-dimensional
structure exhibiting straight interfaces and thus a unique dominant direction
to analyze.
The second geometry we introduce is given by two confocal ellipsoids defining
a coated core in a matrix material. The stiffnesses of the involved materials
are chosen in such a way that the inclusion acts neutrally with respect to a
specific macroscopic strain. This generalizes the (isotropic) 
Hashin structure~\cite{HashinShtrikman}. The curved interfaces make it more difficult
to sample the structure efficiently and the occurring effects can be expected
to be more complex.

\subsection{The laminate structure}\label{sec:laminate}
The probably most simple structure with a predominant direction is a periodic
laminate as shown in Fig.~\ref{fig:laminate:results}
( middle)
.
This structure consists of two isotropic materials alternating in the
direction of lamination
$\vect{n} \in \R^d$ with $\norm{\vect{n}} = 1$ and being constant
perpendicular to it. For this structure Milton~\cite[Section 9.5]{Milton2002}
derives an analytic equation for the effective matrix $\Stiffness^\eff$ that
is given by
\begin{align*}
\left(\mathcal{S}-\mathcal{T}\right)^{-1} &= 
\int_{\T^d} \left(\tilde{\mathcal{S}}(\vect{x})-\mathcal{T}\right)^{-1} \dx \vect{x}
\intertext{with}
\mathcal{S} &\coloneqq \sigma_0 \left( \sigma_0 \Id - \Stiffness^\eff
\right)^{-1},\qquad 
\tilde{\mathcal{S}}(\vect{x}) \coloneqq \sigma_0 \left( \sigma_0 \Id -
\Stiffness(\vect{x}) \right)^{-1},\\
\mathcal{T}_{ijkl} &\coloneqq \frac12 (n_i \delta_{jk} n_l + n_i \delta_{jl} n_k
+n_j \delta_{ik} n_l  + n_j \delta_{il} n_k ) - n_i n_j n_k n_l.
\end{align*}
The choice of the free parameter $\sigma_0$ is explained in 
\cite[Appendix]{KMS:2015Homogenization}. Let $\sigma$ be the largest
eigenvalue of the spectra of the stiffness tensors $\Stiffness(\vect{x})$ then
a choice of $\sigma_0 > \sigma$ ensures that all the inversions necessary to
solve for $\Stiffness^\eff$ can be done.

The resulting strain field $\epsilon$ is, like the structure itself, piecewise
constant and varies only in the direction of lamination. Let the volume
fractions of the two materials be $f_1$ and $f_2$
with $f_1+f_2=1$ and call the corresponding constant strain fields $\epsilon^1$ and
$\epsilon^2$.
The geometry is constant perpendicular to the direction of lamination. Hence
the problem of finding $\epsilon^1$ and $\epsilon^2$ reduces to a
one-dimensional problem.
This gives a system of linear equations involving the macroscopic strain
$\epsilon^0$ to be solved, namely
\begin{equation*}
f_1 \epsilon^1 + f_2 \epsilon^2 = \epsilon^0,\qquad
f_1 \Stiffness^\eff \epsilon^1 + f_2 \Stiffness^\eff \epsilon^2 =
\Stiffness^\eff\epsilon^0.
\end{equation*}

\subsection{The generalized Hashin structure}\label{sec:Hashin}
\begin{figure}\centering
	\begin{subfigure}[b]{.49\textwidth}\centering
	\includegraphics{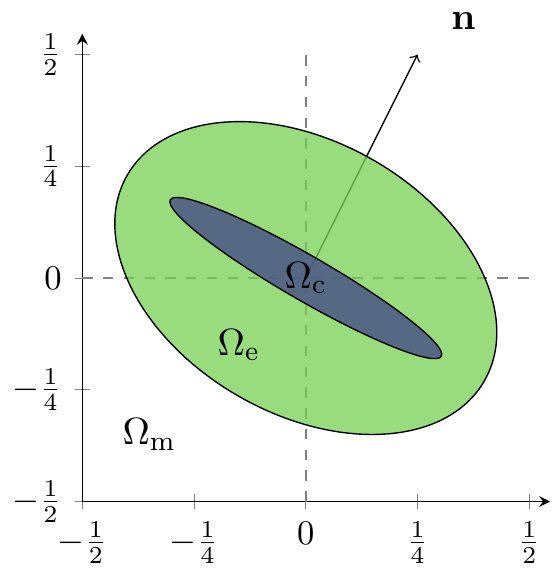}%
	\end{subfigure}
	\begin{subfigure}[b]{.49\textwidth}\centering
			\includegraphics{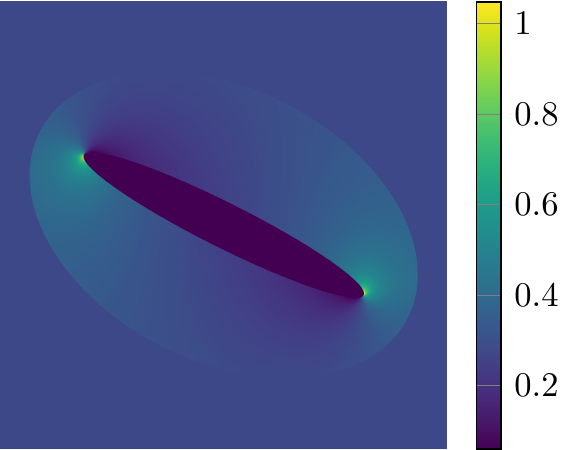}
	\end{subfigure}
	\caption[]{The generalized Hashin structure, \(c_1 = 0.05\), \(c_2 = 0.35\), \(c_3=\infty\), \(\rho_{\mathrm{c}}=0\), \(\rho_{\mathrm{e}}=0.09\), and \(\vect{n} = (\tfrac{1}{2},1)^\tT\) in the \(xy\)-plane, i.e.~the 2D setting is shown in a schematic visualization of the material (left) and the analytic solution \(\epsilon(\vect{x})\) to the elasticity equation (right).}
	\label{fig:Hashin}
\end{figure}

The idea of the generalized Hashin structure due to Hashin and
Shtrikman~\cite{HashinShtrikman} is based on constructing a
inclusion embedded in a matrix material that acts neural to a specific
macroscopic strain, i.e.~the inclusion does not effect the
surrounding stress field. An example of such an inclusion is the assemblage of
coated confocal ellipsoids described by Milton \cite[Section 7.7
ff]{Milton2002} whose derivation we want to follow. A schematic of such a
structure is depicted in Fig.~\ref{fig:Hashin}, left
.

The derivation in Milton is based on an
ellipsoid that is aligned with the coordinate axes. An application of
\cite[Section 8.3]{Milton2002} then allows to generalize this to arbitrary
orientations resulting in the formulae stated in the following. They allow to
predict for a macroscopic strain $\epsilon^0$ that will be given below to
analytically express the resulting strain field $\epsilon$ and the action of
the effective stiffness tensor on this input, $\Stiffness^\eff
\epsilon^0$.

To define the geometry consider confocal ellipsoidal coordinates given for
$\tilde{\vect{x}} \in \R^3$ by the equation
\begin{equation}
\frac{\tilde{x}_1^2}{c_1^2+\rho} + \frac{\tilde{x}_2^2}{c_2^2+\rho} +
\frac{\tilde{x}_3^2}{c_3^2+\rho} = 1
\label{eq:ellipsoidal-coordinates}
\end{equation}
where w.l.o.g. the constants $0 \leq c_1 \leq c_2 \leq c_3 \leq \infty$
determine the relative lengths of the semi-axes of the ellipsoid. A constant
$\rho \geq -c_1^2$ specifies the boundary of an ellipsoid where the lengths of
the semi-axes are given by $l_j(\rho) \coloneqq \sqrt{c_j^2+\rho}$ with $j =
1,2,3$. For a fixed $\tilde{\vect{x}} \in \R^3$ the ellipsoidal radius
$\rho(\tilde{\vect{x}})$ is the uniquely determined largest of the possible three
solutions of \eqref{eq:ellipsoidal-coordinates} and fulfills $\rho(\tilde{x}) \geq
-c_1^2$.

For sake of simplicity, we want to restrict ourselves in this work to prolate
spheroids with $c_1 = c_2 \leq c_3$, oblate spheroids with $c_1 \leq c_2 =
c_3$ and elliptic cylinders. The latter are the limit case $c_3 \rightarrow \infty$.

For these, according to~\cite[Section 7.10]{Milton2002}, the depolarization
factors are given by
\begin{equation*}
d_1(\rho) = d_2(\rho),\qquad
d_2(\rho) = 2-2 d_3(\rho),\qquad
d_3(\rho) = \frac{1-\delta^2}{\delta^2}\Bigl(\frac{1}{2 \delta}
\log\Bigl(\tfrac{1+\delta}{1-\delta}\Bigr) \Bigr),
\end{equation*}
for prolate spheroids, where \(
\delta = \sqrt{1-\tfrac{l_2(\rho)^2}{l_3(\rho)^2}}
\) . Furthermore we have for oblate spheroids using
\( \delta = \sqrt{1-\tfrac{l_1(\rho)^2}{l_2(\rho)^2}}\) the factors
\begin{equation*}
d_1(\rho) = \frac{1}{\delta^2} \Bigl(1-\tfrac{\sqrt{1-\delta^2}}{\delta} \sin^{-1}
\delta \Bigr),\qquad
d_2(\rho) = 2 - 2d_1(\rho),
\qquad
d_3(\rho) =  d_2(\rho),
\end{equation*}
and finally for elliptic cylinders
\begin{equation*}
d_1(\rho) = \frac{l_2(\rho)}{l_1(\rho)+l_2(\rho)},\qquad
d_2(\rho) = \frac{l_1(\rho)}{l_1(\rho)+l_2(\rho)},\qquad
d_3(\rho) = 0.
\end{equation*}

Now let $\rho_{\mathrm{c}}$ and $\rho_{\mathrm{e}}$ be the ellipsoidal radius of the core and the
exterior coating, respectively, cf.~Fig.~\ref{fig:Hashin}%
~(left)%
, with $-c_1^2 < \rho_{\mathrm{c}} < \rho_{\mathrm{e}}$ and with
$l_3(\rho_{\mathrm{e}}) < \frac12$ for $c_3 < \infty$ and $l_2(\rho_{\mathrm{e}}) < \frac12$ for
$c_3 = \infty$, i.e.~the exterior ellipsoid should be contained
in $[-\frac12,\frac12)^3$.

Further, let $\vect{n} \in \R^3$ with $\norm{\vect{n}} = 1$ be the direction the shortest
semi-axis of the ellipsoid with length $l_1(\rho)$. Define the rotation matrix
that transforms the vector $(1,0,0)^\tT$ to $\vect{n}$ by
\begin{equation*}
\mat{R} \coloneqq \begin{pmatrix} 1 & -n_2 & -n_3 \\ n_2 & 1 & 0\\ n_3 & 0 & 1
\end{pmatrix}+ \frac{1-n_1}{\sqrt{n_2^2+n_3^2}}
\begin{pmatrix} -n_2^2-n_3^2 & 0 & 0\\ 0 & -n_2^2 & -n_2n_3\\0 & -n_2n_3 &
-n_3^2 \end{pmatrix}.
\end{equation*}

Then the core, the exterior coating and the surrounding matrix are given by
\begin{align*}
\Omega_{\mathrm{c}} &\coloneqq \left\lbrace \vect{x} \in \R^3 : \rho(\mat{R}^{-1}\vect{x}) \leq
\rho_{\mathrm{c}} \right\rbrace,\\
\Omega_{\mathrm{e}} &\coloneqq \left\lbrace \vect{x} \in \R^3 : \rho_{\mathrm{c}} < \rho(\mat{R}^{-1}\vect{x}) \leq
\rho_{\mathrm{e}} \right\rbrace,\\
\Omega_{\mathrm{m}} &\coloneqq \bigl[ -\tfrac{1}{2}, \tfrac{1}{2} \bigr)^3 \setminus \left(
\Omega_{\mathrm{c}} \cup \Omega_{\mathrm{e}} \right),
\end{align*}
respectively. With
\begin{equation*}
f(\rho) \coloneqq \frac{\sqrt{g(\rho_{\mathrm{c}})}}{\sqrt{g(\rho)}}
\quad\text{ and }\quad
g(\rho) \coloneqq \begin{cases} (c_1^2+\rho)(c_2^2+\rho)(c_3^2+\rho),
&\text{for } c_3 < \infty,\\
(c_1^2+\rho)(c_2^2+\rho),
&\text{else}, \end{cases}
\end{equation*}
the volume fraction of $\Omega_{\mathrm{c}}$ in the coated ellipsoid is
given by $f(\rho_{\mathrm{e}})$.

We assume that the material in the core and in the exterior coating behave
isotropically, i.e.~they are described by stiffness matrices of the form 
$\Stiffness_{ijkl} = \lambda \delta_{ij} \delta_{kl} + \mu(\delta_{ik}
\delta_{jl}+\delta_{il} \delta_{jk})$. The parameter $\lambda$ is Lam\'e's
first parameter and $\mu$ is the shear modulus. We denote the parameters in
the core by $\lambda_{\mathrm{c}}$ and $\mu_{\mathrm{c}}$ and in the exterior coating by $\lambda_{\mathrm{e}}$
and $\mu_{\mathrm{e}}$, respectively.
Further, the bulk modulus of an isotropic material is given as $\kappa
\coloneqq \lambda + \frac23 \mu$.

Following \cite[Section 7.9]{Milton2002} we impose a macroscopic strain of
\begin{equation*}
\epsilon^0 = \mat{R}\left(\frac{3 \kappa_{\mathrm{e}} +4 \mu_{\mathrm{e}}}{9(\kappa_{\mathrm{c}}-\kappa_{\mathrm{e}})} \Id +
(1-f(\rho_{\mathrm{e}}))\mat{S}(\rho_{\mathrm{e}})\right) \mat{R}^\tT
\end{equation*}
with
\begin{align*}
\mat{S}(\rho) &\coloneqq (1-f(\rho))^{-1} (\mat{D}(\rho_{\mathrm{c}}) - f(\rho) \mat{D}(\rho))\quad\text{ and }\quad
\mat{D}(\rho) \coloneqq \diag(d_i(\rho)_{i=1,2,3}).
\end{align*}

The effective matrix of the structure and likewise ---to ensure the
neutrality of the inclusion--- the stiffness matrix of the
matrix material $\Omega_{\mathrm{m}}$ are given by their action on the macroscopic strain
as 
\begin{equation*}
\Stiffness^\eff \epsilon^0 \coloneqq
\mat{R}\Bigl(  	\tfrac{\kappa_{\mathrm{e}}}{\kappa_{\mathrm{c}}-\kappa_{\mathrm{e}}}(\kappa_{\mathrm{c}}+\frac43 \mu_{\mathrm{e}}) 
+ \frac43 \mu_{\mathrm{e}} f(\rho_{\mathrm{e}}) \Bigr) \Id  \mat{R}^\tT
 + \mat{R}\frac23 \mu_{\mathrm{e}} (1-f(\rho_{\mathrm{e}}))(3 \mat{S}(\rho_{\mathrm{e}}) - \Id) \mat{R}^\tT.
\end{equation*}
The resulting strain field is then for $\vect{x} \in \Omega_{\mathrm{c}}$ given as
\begin{equation}
\epsilon(\vect{x}) = \mat{R}\frac{3 \kappa_{\mathrm{e}} + 4 \mu_{\mathrm{e}}}{9(\kappa_{\mathrm{c}}-\kappa_{\mathrm{e}})} \Id  \mat{R}^\tT,
\end{equation}
and for $\vect{x} \in \Omega_{\mathrm{m}}$ as \(\epsilon(\vect{x}) = \epsilon^0\), respectively. The strain field is constant in the core and in the matrix material.
In the external coating we have for $\vect{x} \in \Omega_{\mathrm{e}}$
\begin{align*}
\epsilon(\vect{x})
	&= \mat{R}
	\biggl(
		\tfrac{3 \kappa_{\mathrm{e}} + 4 \mu_{\mathrm{e}}}{9(\kappa_{\mathrm{c}}-\kappa_{\mathrm{e}})}
\Id
+
\mat{D}(\rho_{\mathrm{c}}) - f(\rho(\tilde{\vect{x}}))
\mat{D}(\rho(\tilde{\vect{x}}))
+
\tfrac{\sqrt{g(\rho_{\mathrm{c}})}}{2} \vect{q}(\tilde{\vect{x}})
\Grad_{\tilde{\vect{x}}}^T \rho(\tilde{\vect{x}})
	 \biggr)\mat{R}^\tT,
\intertext{with $\tilde{\vect{x}} \coloneqq \mat{R}^{-1} \vect{x}$ and }
\vect{q}(\tilde{\vect{x}})_i &=  \frac{\tilde{x}_i}{(c_i^2+\rho(\tilde{\vect{x}}))
\sqrt{g(\rho(\tilde{\vect{x}}))}},
\quad
\left(\Grad_{\tilde{\vect{x}}} \rho(\tilde{\vect{x}}) \right)_i = \frac{2
\tilde{x}_i}{c_i^2+\rho(\tilde{\vect{x}})} \Biggl(\sum_{j=1}^3
\frac{\tilde{x}_i}{c_i^2+\rho(\tilde{\vect{x}})}\Biggr)^{-1},
\intertext{ $i=1,2,3$. For the case $c_3 = \infty$ we additionally have $\vect{q}(\tilde{\vect{x}})_3 = 0$,}
 \left(\Grad_{\tilde{\vect{x}}} \rho(\tilde{\vect{x}}) \right)_i
 &=
  \frac{2\tilde{x}_i}{c_i^2+\rho(\tilde{\vect{x}})} \Biggl(\sum_{j=1}^2
 \frac{\tilde{x}_i}{c_i^2+\rho(\tilde{\vect{x}})}
	\Biggr)^{-1},\quad i=1,2,\text{ and }\left(\Grad_{\tilde{\vect{x}}} \rho(\tilde{\vect{x}}) \right)_3 = 0.
\end{align*}

\begin{example}\label{ex:hashin}
Let an  ellipse in 2D be given, i.e.~we set the third dimension to be constant by \(c_3=\infty\), by choosing \(c_1 = 0.05\), \(c_2 = 0.35\). Choose further \(\rho_c = 0\) and \(\rho_e = 0.09\) and introduce a rotation of \(60^\circ\) counter clockwise by setting \(\vect{n} = (\tfrac{1}{2},1,0)^\tT\). The resulting geometry is shown in ~\ref{fig:Hashin}%
~(left)%
, the analytic solution \(\epsilon(\vect{x})\) to the elasticity equation in%
~\ref{fig:Hashin} (right).
\end{example}
\section{Numerics}\label{sec:numerics}

The algorithms used in this paper are implemented in \lstinline!MatLab! R2015b
in a modular and fast way using vectorization. For the Fourier transform on
arbitrary patterns we employ the multivariate periodic anisotropic wavelet library
(MPAWL)\cite{Bergmann2013MPAWL}, which was recently ported to
Matlab\footnote{see~\url{https://github.com/kellertuer/MPAWL-Matlab}}. It uses
Matlab's internal~\lstinline!fftn! command to apply the fast Fourier transform.
All tests were run on a MacBook Pro running Mac OS X 10.11.5, Core i5, 2.6 GHz,
with 8 GB RAM using \textsc{Matlab} 2016a and the clang-700.1.76 compiler.

\subsection{Hashin}

Consider the geometry of the coated ellipsoid as described in
Section~\ref{sec:Hashin} with the parameters in Example~\ref{ex:hashin} seen as
a two-dimensional problem, i.e.~sampled with only one point in
$x_3$-direction.
This structure is strongly orthotropic with the dominant directions being
$\vect{n}$ and $\vect{n}^\perp$. To analyze the influence of the sampling
matrix $\mat{M}$ we compare the relative $\ell^2$-error of the strain $\epsilon$
compared with the analytic solution and the relative
error in the effective stiffness tensor sampled on \(\Pattern(\mat{M})\), i.e we define
\begin{equation*}
e_{\ell^2}(\mat{M}) \coloneqq
\frac{\norm{\epsilon-\tilde{\epsilon}}[l^2(\Pattern(\mat{M}))]}
{\norm{\tilde{\epsilon}}[l^2(\Pattern(\mat{M}))]}
\quad\text{ and }\quad
e_{\eff}(\mat{M}) \coloneqq
\frac{\norm{\Stiffness^\eff \epsilon^0 - \tilde{\Stiffness}^\eff
\epsilon^0}[2]}{\norm{\tilde{\Stiffness}^\eff \epsilon^0}[2]},
\end{equation*}
where $\epsilon$ and $\Stiffness^\eff \epsilon^0$ are the numerical solutions
obtained by Algorithm~\ref{alg:fixed-point}, and $\tilde{\epsilon}$ and
$\tilde{\Stiffness}^0 \epsilon^0$ are the analytic solutions.

The pattern matrices are parametrized by
\begin{equation}
\mat{M}_{j,k,\alpha} \coloneqq
\begin{pmatrix}
2^j & \alpha k & 0\\
(1-\alpha)k & 2^{14-j} & 0\\
0 & 0 & 1
\end{pmatrix}
\label{eq:M-shear}
\end{equation}
with $j \in \lbrace 7,\dots,9 \rbrace$, $\alpha \in \lbrace 0,1 \rbrace$ and $k
\in 16 \cdot \lbrace -32, \dots, 32 \rbrace$. For all parameters these matrices
have determinant $2^{14}$, i.e.~the number of sampling points stays constant. The parameter $k$ shears the pattern, $\alpha$ determines the direction of the
shearing and $j$ controls the refinement in the direction of the pattern basis
vectors~\(\vect{y}_j\), \(j=1,\ldots,d_{\mat{M}}\), cf.~\eqref{eq:patternBasisV}. This induces both an anisotropy or preference of direction in the pattern~\(\Pattern(\mat{M})\) as well as for the basis vectors~\(\vect{h}_j\) of the corresponding generating set~\(\generatingSet(\mat{M}^\tT)\), cf.~\eqref{eq:genSetBasisV}, representing the frequencies.

Further, we use pattern matrices of the form
\begin{equation*}
\tilde{\mat{M}}_j \coloneqq \frac12 
\begin{pmatrix} 1&1&0\\-1&1&0\\0&0&2\end{pmatrix}
\begin{pmatrix}2^j&0&0\\0&2^{14-j}&0\\0&0&1\end{pmatrix}
\end{equation*}
corresponding to a rotation of the grid in the direction of $\vect{n}$ with
$\det(\tilde{\mat{M}}_j) = 2^{13}$.

\begin{figure}\centering
\includegraphics{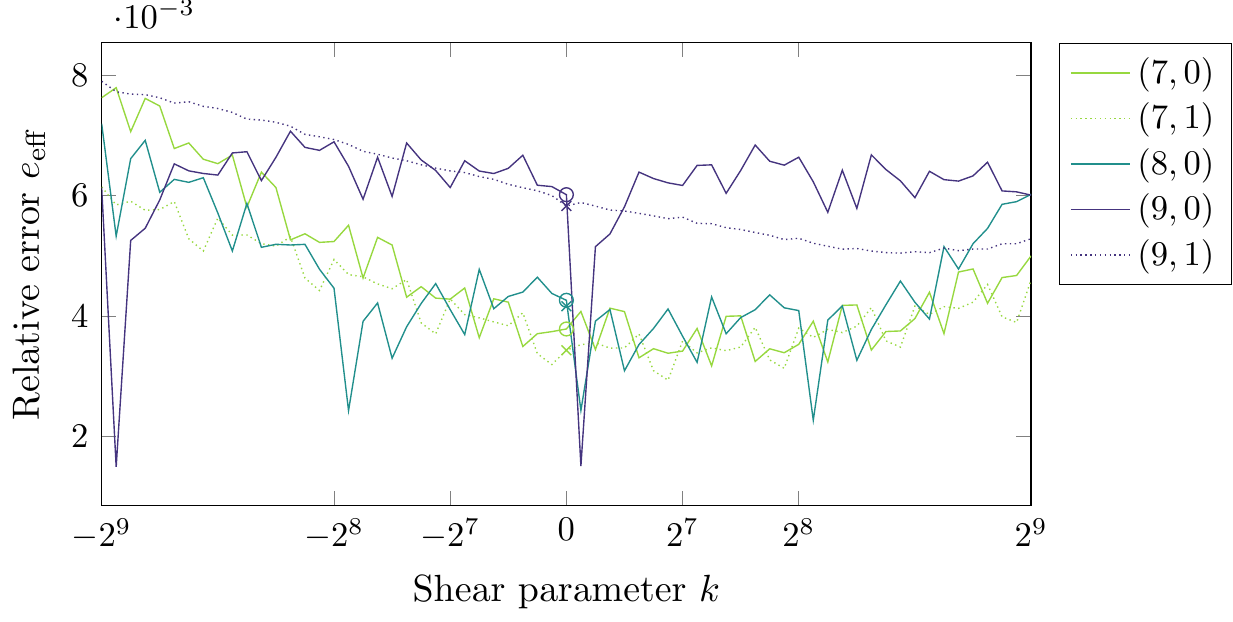}
\caption{Shear parameter against relative error in the effective stiffness
tensor. The corresponding matrices for the curve $(j,\alpha)$ at point $k$ is
$\mat{M}_{j,k,\alpha}$ from~\protect\eqref{eq:M-shear}. Circles mark the
results for matrices $\mat{M}_{j,0,0}$ and crosses stand for matrices of the
form $\tilde{\mat{M}}_j$.}
\label{fig:ellipsoid2Drel}
\end{figure}
In Figure~\ref{fig:ellipsoid2Drel} the effect of the shearing parameter $k$ on
$e_{\eff}$ is depicted, neglecting curves that do not perform better than
$\mat{M}_{7,k,\alpha}$ with $\mat{M}_{9,k,1}$ as an example for such a curve.
For each value of~\(j\)
we choose one color (brightness) and indicate \(\alpha=0\) by a solid,
\(\alpha=1\) by a dotted line, and indicate \(\tilde{\mat{M}}_j\), \(j=7,8,9\), by
crosses as well as the corresponding diagonal
matrices~\(\mat{M}_{j,0,0}=\mat{M}_{j,0,1}\) by circles.
The circles correspond to the classical rectangular “pixel sampling grid”.

The reference point for the analysis of the results is the unsheared
matrix $\mat{M}_{7,0,0}$, i.e.~the standard equidistant Cartesian grid, giving
an error of $3.8 \cdot 10^{-3}$.
Shearing this matrix in either direction results in a larger errors, e.g.~ with
twice the error for $\mat{M}_{7,-512,0}$ and an error of $5 \cdot 10^{-3}$ for
$\mat{M}_{7,512,0}$.

In contrast $\mat{M}_{8,k,0}$ behaves similarly to $\mat{M}_{7,k,0}$ with the
exception of $k_{-1} = -2^8+16$, $k_{0} = 16$ and $k_{1} = 2^8+16$, where the
error is $2.4 \cdot 10^{-3}$ for the first two and $2.2 \cdot 10^{-3}$ for the
third point. Note that the matrices possess the same pattern due to
\(\mat{M}_{8,k_{-1},0} \sim_{\Pattern} \mat{M}_{8,k_0,0} \sim_{\Pattern}
\mat{M}_{8,k_1,0}\).

\begin{figure*}
		\begin{subfigure}{\textwidth}\centering
			\includegraphics{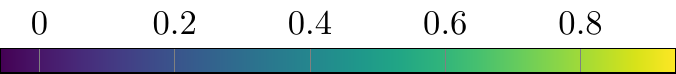}
		\end{subfigure}
	\\[.25\baselineskip]
	\begin{subfigure}[t]{.24\textwidth}\centering
		\includegraphics[width=0.95\textwidth]{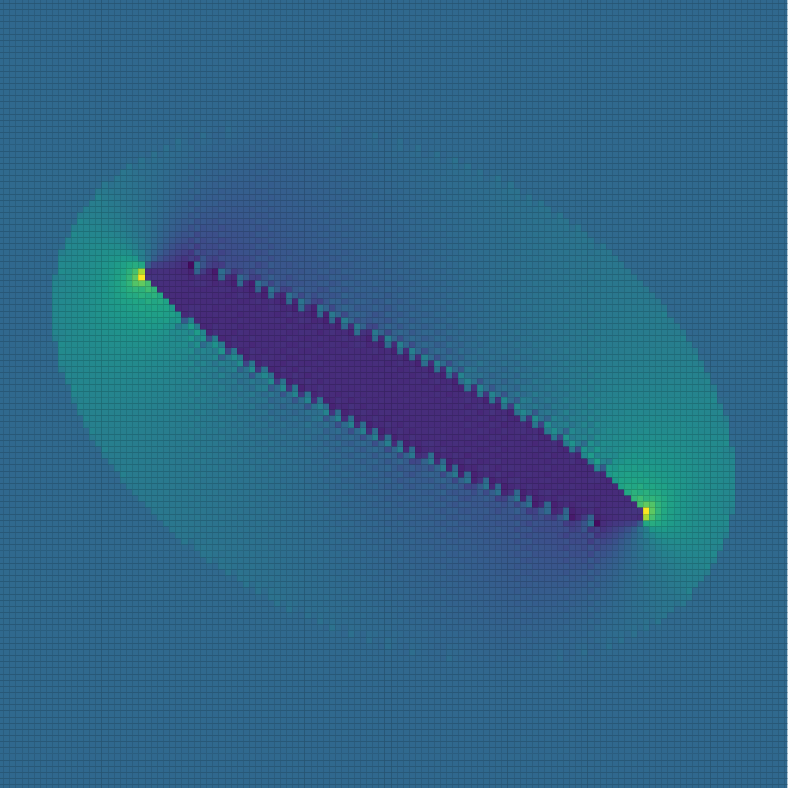}
	\end{subfigure}
	\begin{subfigure}[t]{.24\textwidth}\centering
		\includegraphics[width=0.95\textwidth]{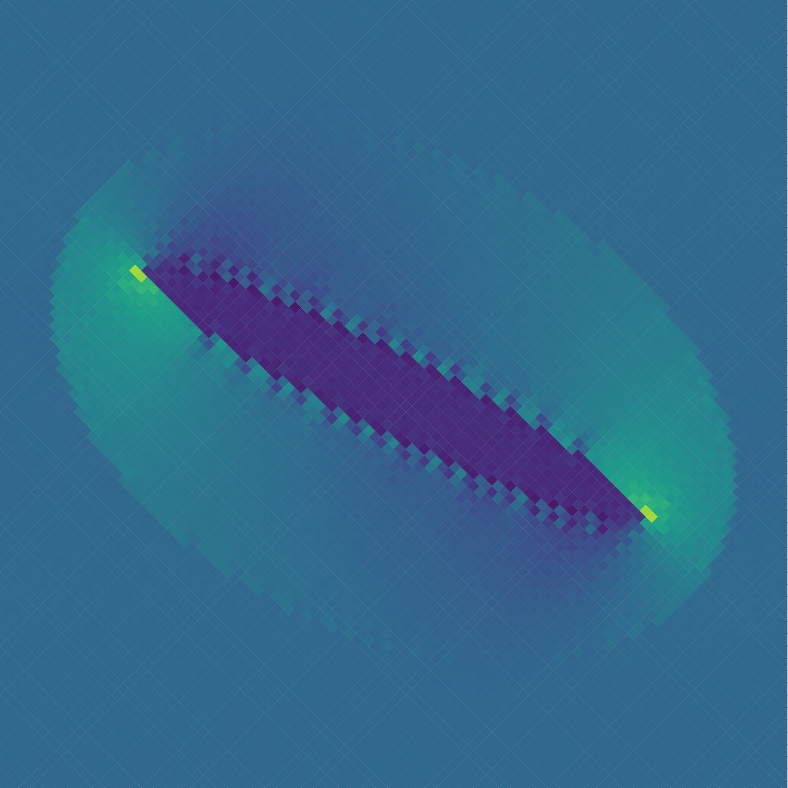}
	\end{subfigure}
	\begin{subfigure}[t]{.24\textwidth}\centering
		\includegraphics[width=0.95\textwidth]{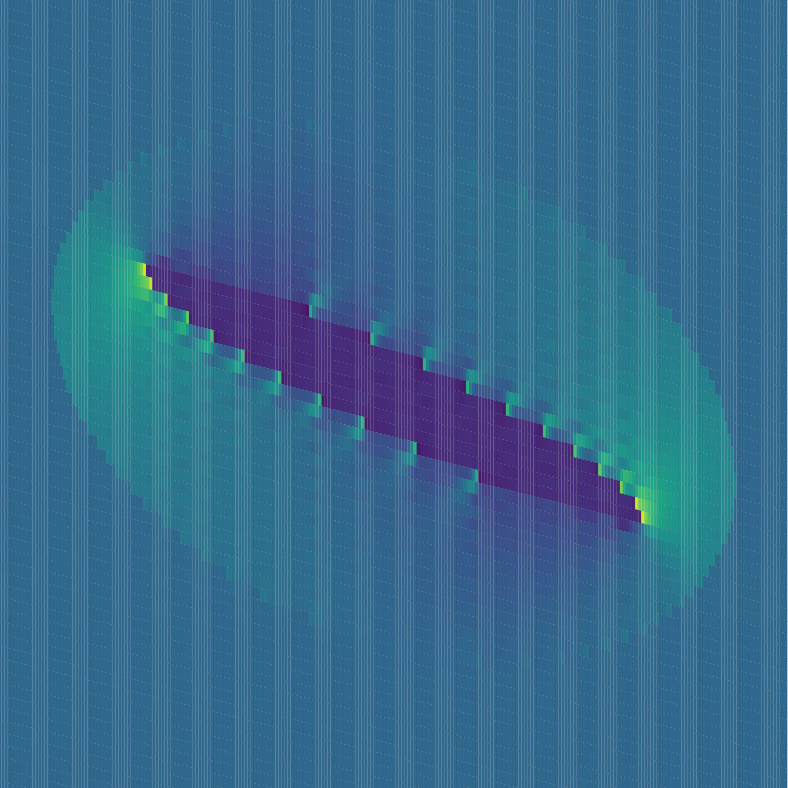}
	\end{subfigure}
	\begin{subfigure}[t]{.24\textwidth}\centering
		\includegraphics[width=0.95\textwidth]{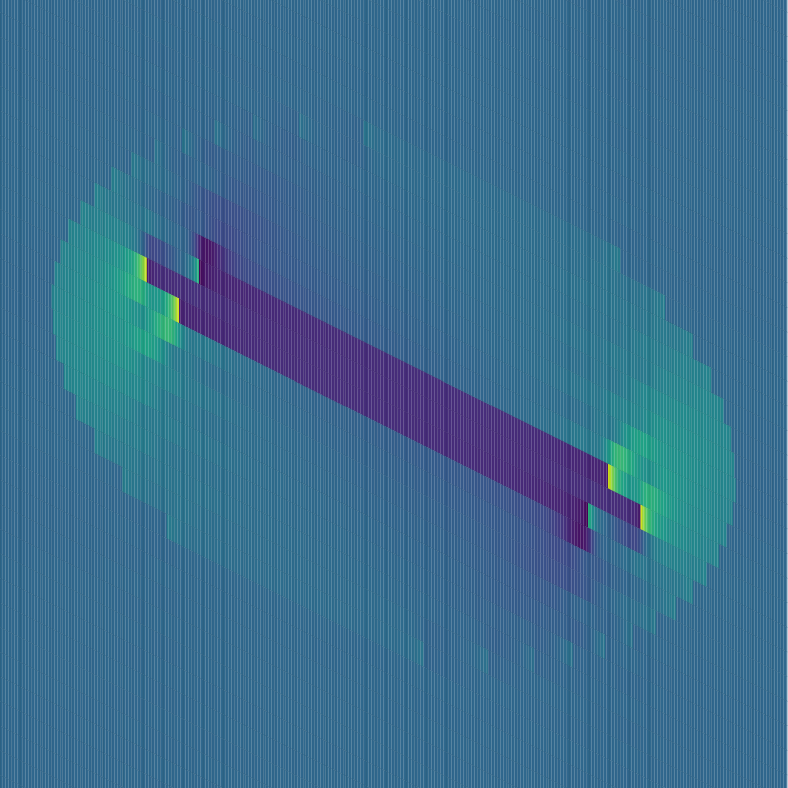}
	\end{subfigure}
	\\
	\begin{subfigure}[t]{.24\textwidth}\centering
		\includegraphics[width=0.95\textwidth]{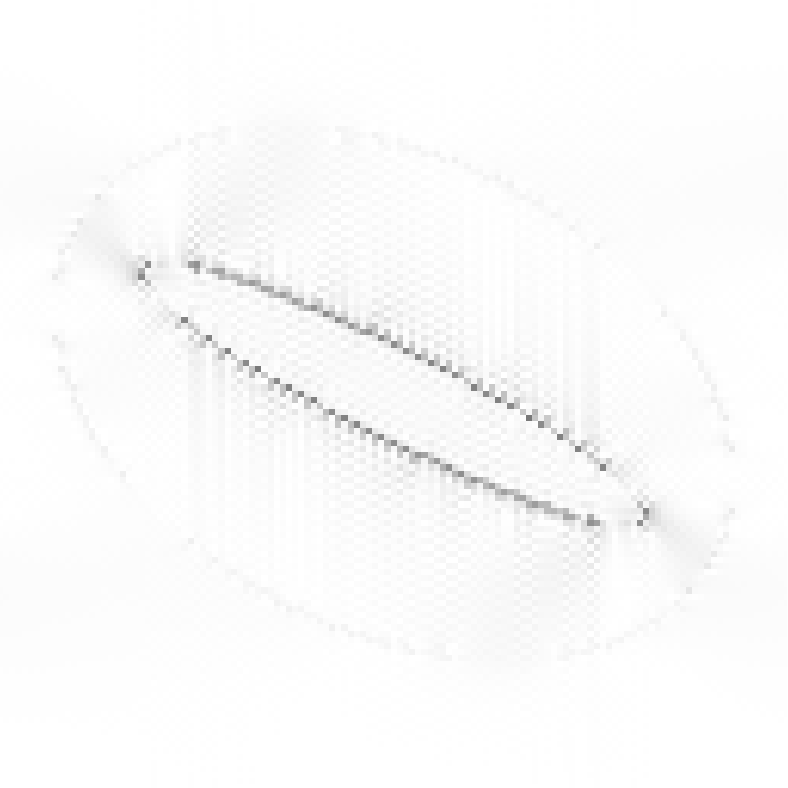}
	\end{subfigure}
	\begin{subfigure}[t]{.24\textwidth}\centering
		\includegraphics[width=0.95\textwidth]{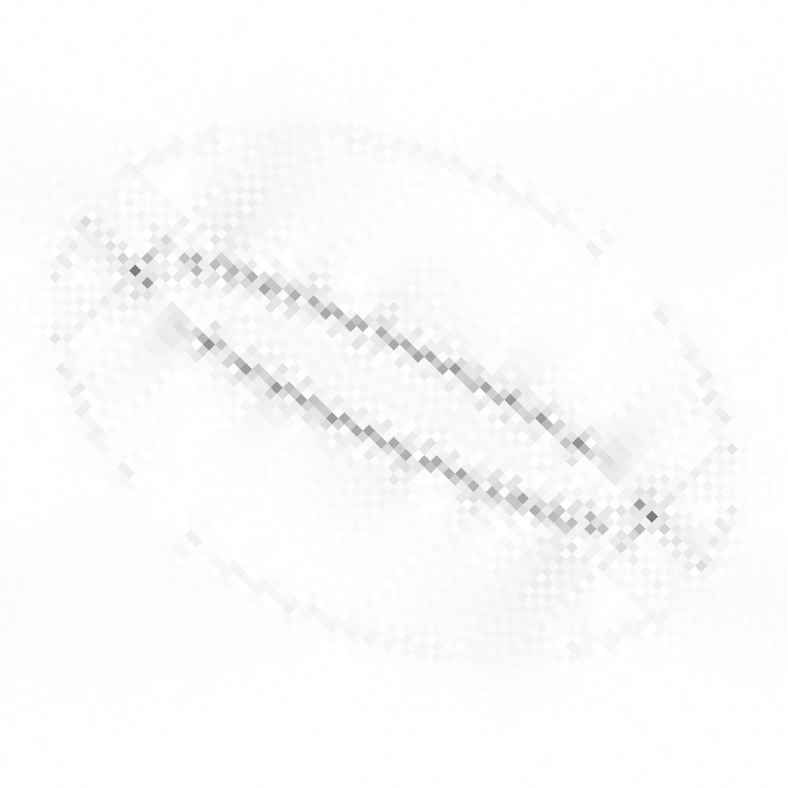}
	\end{subfigure}
	\begin{subfigure}[t]{.24\textwidth}\centering
		\includegraphics[width=0.95\textwidth]{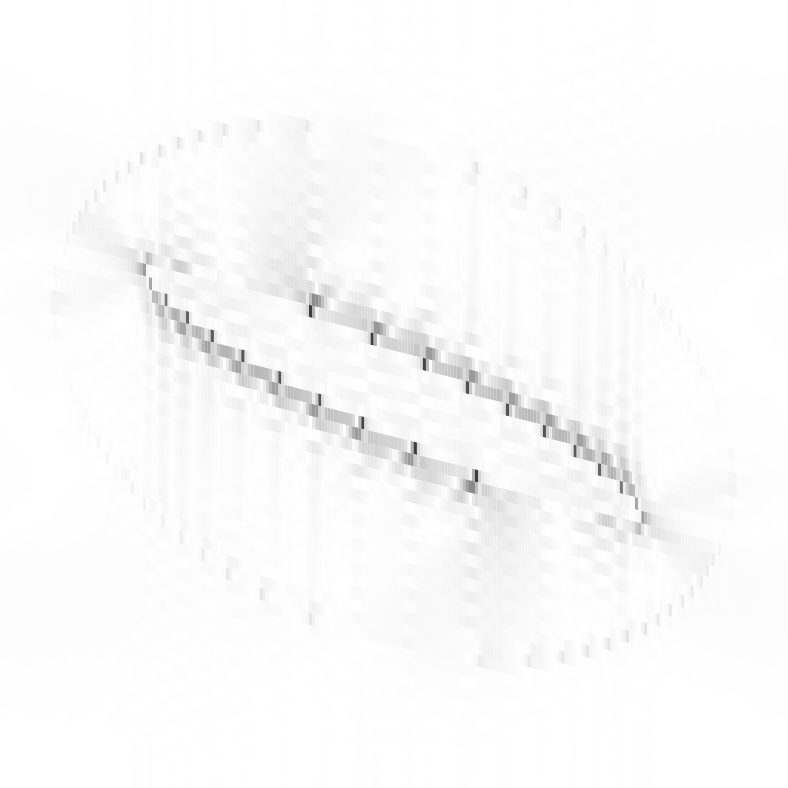}
	\end{subfigure}
	\begin{subfigure}[t]{.24\textwidth}\centering
		\includegraphics[width=0.95\textwidth]{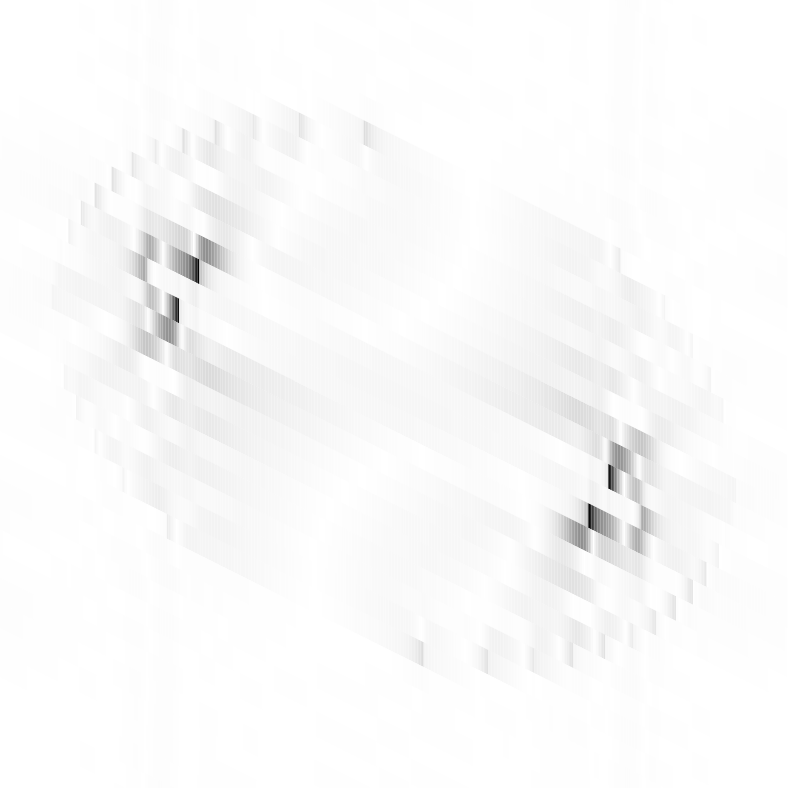}
	\end{subfigure}
	\\[.25\baselineskip]
		\begin{subfigure}{\textwidth}\centering
			\includegraphics{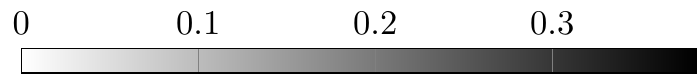}
		\end{subfigure}
	\caption[]{
	The strain field \(\epsilon_{11}\) (first row) using the colormap above and
 	its corresponding \(e_{\log}\)-error (second row) given by
 	\(e_{\log}=\log(1+\lvert\epsilon-\tilde\epsilon\rvert)\) using the colormap
 	below. The matrices inducing the sampling pattern points
 	are (left to right):
	\(\mat{M}_{7,0,0}=\bigl(\begin{smallmatrix}2^7&0\\0&2^7\end{smallmatrix}\bigr)\),
	\(\tilde{\mat{M}}_{7}=\bigl(\begin{smallmatrix}2^6&2^6\\-2^6&2^6\end{smallmatrix}\bigr)\),
 	\(\mat{M}_{8,16,0}
 	=\bigl(\begin{smallmatrix}2^8&0\\2^4&2^6\end{smallmatrix}\bigr)
 	\), and
	\(\mat{M}_{9,16,0}
	=\bigl(\begin{smallmatrix}2^9&0\\2^4&2^5\end{smallmatrix}\bigr)
	\)%
.}\label{fig:Hashin:Ex}
\end{figure*}

This effect is even more dominant for $\mat{M}_{9,k,0}$ that gives an error of
around $6 \cdot 10^{-3}$ almost everywhere except for $k = -2^9+16$ and $k=16$
where it drops to $1.5 \cdot 10^{-3}$. These matrices are also congruent with
respect to pattern congruence \(\sim_{\Pattern}\).

The strain field and the pointwise error for the $\epsilon_{11}$-component of
$\mat{M}_{7,0,0}$ are depicted in
Fig.~\ref{fig:Hashin:Ex}%
, first column. As the
ellipsoid is not aligned with the pattern, the interface along the long side of
the ellipsoid is not resolved well. As the sine functions are not perpendicular
to this interface the Gibbs phenomenon dominates the $\ell^2$-error. Likewise the
strain peak at the tip of the ellipsoid is not captured correctly and worsens
$e_\eff$ tremendously.

If we choose a matrix with a certain shear, e.g.~\(\mat{M}_{8,16,0}\), where the
strain \(\epsilon_{11}\) and the corresponding error are shown in
Fig.~\ref{fig:Hashin:Ex}
 (third column)
, we obtain a quite small error in the
effective stiffness tensor and resolve the strain peaks correctly. As the
pattern, however, is not aligned with the ellipsoid, the interfaces show large
errors and inwards and outwards facing corners result in a very large
$\ell^2$-error.

For the matrices $\mat{M}_{9,-2^9+16,0}$ and $\mat{M}_{9,16,0}$ the
pattern is sheared in such a way that it is aligned with the
ellipsoid, c.f. Fig.~\ref{fig:Hashin:Ex}%
, last column,
refined in along the longer semi-axis.
The strain field $\epsilon_{11}$ is characterized by 
slowly changing values in the direction of the shorter semi-axis and especially
rapidly changing and high strains at the tips of the inner ellipsoid. Therefore,
to get small errors, we need only few points in the direction of the shorter
semi-axis, in this case only \(32\) points.
The high strains at the tips of the ellipsoid, however, require a high
resolution like $512$ sampling points in this case. This leads to a lower
approximation of the edges orthogonal to the smaller semi-axis and hence the
Gibbs phenomenon increases the $\ell^2$-norm, cf.~Fig.~\ref{fig:Hashin:Ex}%
, first column,
for the log-error.
This effect concentrates around the tips of the ellipsoid resulting in a core almost shaped like a rectangle. The averaging done to
compute the effective stiffness tensor cancels these errors and does therefore
not influence the error.

\begin{table}\centering
\begin{tabular}{lcc}\toprule
\(\mat{M}\) & $e_\eff(\mat{M})$ & $e_{\ell^2}(\mat{M})$\\
\midrule
$\mat{M}_{7,0,0}$       & 0.0038 & 0.043\\
$\tilde{\mat{M}}_7$     & 0.0034 & 0.047\\\midrule
$\mat{M}_{8,-2^8+16,0}$ & 0.0024 & 0.054\\
$\mat{M}_{8,16,0}$      & 0.0024 & 0.050\\
$\mat{M}_{8,2^8+16,0}$  & 0.0023 & 0.048\\\midrule
$\mat{M}_{9,-2^9+16,0}$ & \textbf{0.0015} & 0.054\\
$\mat{M}_{9,15,0}$      & 0.0019 & 0.052\\
$\mat{M}_{9,16,0}$      & \textbf{0.0015} & 0.054\\
$\mat{M}_{9,17,0}$      & 0.0024 & 0.054\\
	\midrule
$\mat{M}_{7,2^8,1}$     & 0.0038 & 0.023\\
$\mat{M}_{7,2^8+16,1}$  & 0.0036 & \textbf{0.022}\\\bottomrule
\end{tabular}
\caption{Relative effective stiffness and \(\ell^2\)-errors for several shearing matrices. While matching the direction, \(\mat{M}_{9,16,0}\) reduces the \(e_{\text{eff}}\) error, the \(e_{\ell^2}\)-error is reduced tremendously by e.g.~shearing the standard grid. This can be seen by looking at $\mat{M}_{7,2^8+16,1}$.}
\label{tab:ellipsoid}
\end{table}

Table~\ref{tab:ellipsoid} shows both the \(\ell^2\)- and the \(e_{\eff}\)-errors for several matrices. The smallest value
of $e_\eff$ is reached for $\mat{M}_{9,-2^9+16,0}$ and $\mat{M}_{9,16,0}$.
The error for shear parameters around $\mat{M}_{9,16,0}$ get slightly worse by
changing to $k=15$ or $k=17$,
respectively. The smallest $\ell^2$-errors can be obtained by $\mat{M}_{7,2^8,1}$ and
$\mat{M}_{7,2^8+16,1}$, giving errors only half as large as in the standard
case. The pattern $\mat{M}_{7,2^8,1}$ involves taking a standard pixel grid for
sampling and then shearing the unit cell in such a way that it is aligned with
the ellipsoid. This gives a balance between resolving the interfaces and
capturing the strain peak. The alignment of the sine ansatz functions with the
ellipsoid reduces the Gibbs phenomenon and the comparatively large number of
frequencies used in the direction of the strain peaks allows for small errors
there. The matrix $\mat{M}_{7,2^8+16,1}$ further shears the pattern by a small
amount and resolves the strain peak better while preserving the good
resolution of the ellipsoid.

Subsampling of the complete rectangular grid of one \(\mat{M}_{j,0,0}\) on the
so-called quincunx pattern induced by the matrix \(\tilde{\mat{M}}_j\) shown as circles for the former and crosses for the latter matrices in
Fig.~\ref{fig:ellipsoid2Drel}, gives slightly smaller values of $e_\eff$.
The $\ell^2$-errors like for $\tilde{\mat{M}}_7$, cf.
Figs.~\ref{fig:Hashin:Ex}%
, second column,
also decrease even especially when taking into account that only half the sampling values are used.

\subsection{Subsampling}
	\begin{figure*}%
		\centering
		\begin{subfigure}{\textwidth}\centering
			\includegraphics{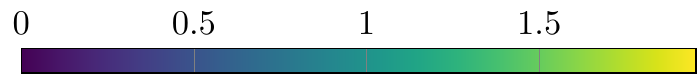}
		\end{subfigure}
	\\[.25\baselineskip]
		\def\axisdim{.975\textwidth}
		\begin{subfigure}{.32\textwidth}\centering
			\includegraphics{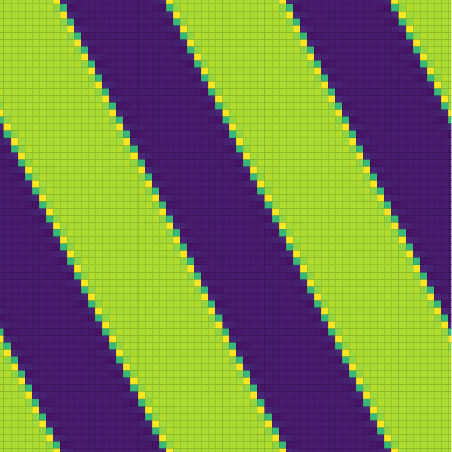}
		\end{subfigure}
		\begin{subfigure}{.32\textwidth}\centering
			\includegraphics{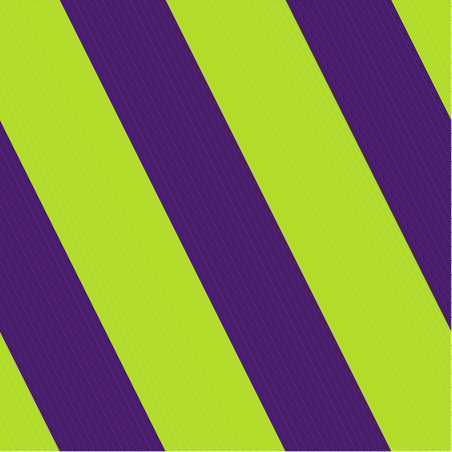}
		\end{subfigure}
		\begin{subfigure}{.32\textwidth}\centering
			\includegraphics{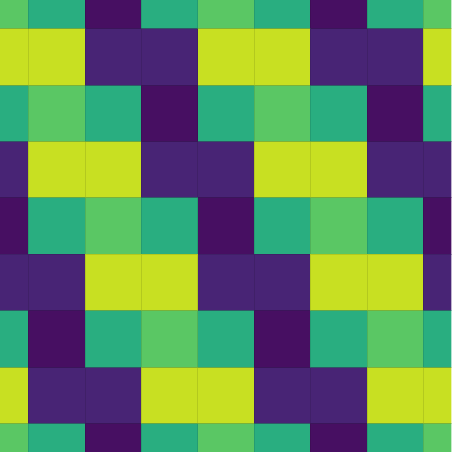}
		\end{subfigure}
		\caption{A comparison of 3 different sampling sets of the laminate material: a full \(64\times 64\) (left, \(e_{\ell^2}(\mat{M}_a)=0.109\), 94 Iterations), an adapted rank-1-lattice having \(64\) points (middle, \(e_{\ell^2}(\mat{M}_b)=0.0348\), 9 It.) and the comparison pixel grid also having \(64\) points, i.e. \(8\times 8\) pixel (right, \(e_{\ell^2}(\mat{M}_c)=0.2909\), 89 It.). The latter two are both subsamplings of the first.}
		\label{fig:laminate:results}
	\end{figure*}%

We study possibilities to subsample given (large) data, when certain
directional information is given, e.g.~the laminate from Section~\ref{sec:laminate} and a given normal vector of \(\vect{n} = (1,\tfrac{1}{2})^\tT\). We choose a matrix \(\mat{N}\) being a factorization of a pixel grid and of the form
\[
	\mat{M}
	=
	\begin{pmatrix}
		a&0\\0&a
	\end{pmatrix}
	=
	\begin{pmatrix}
		1&-\frac{a}{2}\\0&a
	\end{pmatrix}
	\begin{pmatrix}
		a & \frac{a}{2}\\0&1
	\end{pmatrix}
	\eqqcolon \mat{J}\mat{N}
\]
and hence \(\Pattern(\mat{N}) \subset \Pattern(\mat{M})\) if \(\frac{a}{2}\) is an integer.
Note that in the sub pattern the most dominant direction in the Fourier domain
is the direction orthogonal to the edge direction.
Setting \(a=2^{6}\) the pixel sampling contains just \(2^{12} = 4096\) data
items, a size where the dominant numerical effects like the Gibbs phenomenon can
still be observed. We compare the following patterns: first we sample on the full grid,
i.e.~we choose \(\mat{M}_{\text{a}} = \bigl(\begin{smallmatrix}64&0\\0&64\end{smallmatrix}\bigr)\). We consider a directional sub lattice given by~\(\mat{M}_\text{b} = 
\bigl(\begin{smallmatrix}
	64&32\\0&1
\end{smallmatrix}\bigr)
\) following the construction above and resulting in $64$ data points, only the square root of the number of points of the full grid.
By construction we have \(\Pattern(\mat{M}_\text{b}) \subset\Pattern(\mat{M}_\text{a})\) and that the spanning vector~\(\vect{y}_1\) of this rank-\(1\)-lattice, is orthogonal to the edges present in the material.
This can for example be done by examining the large(r) dataset with respect to its edge directions and subsampling accordingly.
We study the effect of these in contrast to \(\mat{M}_\text{c} = \bigl(
\begin{smallmatrix}
	8&0\\0&8
\end{smallmatrix}
\bigr)\) having also \(64\) points. Clearly this is also a sub lattice of \(\mat{M}_\text{a}\). We employ Algorithm~\ref{alg:fixed-point} and a Chauchy criterion for stopping, i.e.~\(\tfrac{\norm{\epsilon^{(n+1)}-\epsilon^{(n)}}}{\norm{\epsilon^{(0)}}}\)
with a threshold of \(10^{-9}\). The results are shown in
Figure~\ref{fig:laminate:results}%
.
	
The full grid based approach of~\(\mat{M}_\text{a}\) shown
in Fig.~\ref{fig:laminate:results}%
~(left)
suffers from the well known Gibbs phenomenon. 
The subsampling on~\(\Pattern(\mat{M}_\text{b})\) in Fig~\ref{fig:laminate:results} (middle)
reduces the number of iterations tremendously from \(94\) 
to only \(9\) and the computational times from \(41.2\) to only \(4.5\) seconds. This is not the case for the pixel grid subsampling given by~\(\Pattern(\mat{M}_\text{c})\) shown in Fig~\ref{fig:laminate:results} (right)
which still requires \(89\) iterations and \(46.1\) seconds.

Looking at the \(\ell^2\)-error depicted in the captions of Fig.~\ref{fig:laminate:results}, reducing
the number of points from \(64^2\) in \(\Pattern(\mat{M}_a)\) to \(64\) in \(\Pattern(\mat{M}_b)\) \emph{also} reduces the error by a factor of roughly 3.1. The tensor pixel grid of \(\mat{M}_{\text{c}}\) with \(8\times 8 = 64\) is by a factor of \(8\) worse than our new anisotropic approach given by~\(\mat{M}_b\). 

Finally, we analyse the error of the effective stiffness \(e_{\text{eff}}(\mat{M}_x)\), \(x\in\{a,b,c\}\).
The values are \(e_{\text{eff}}(\mat{M}_a)=0.0042\), \(e_{\text{eff}}(\mat{M}_b)=0.0134\), and \(e_{\text{eff}}(\mat{M}_c)=0.0495\). Hence having only \(64\) times the number of points as \(\Pattern(\mat{M}_{\text{b}})\) the first example is only about a factor of~\(3.2\) times better.
This pattern yields an error that is by a factor of roughly \(3.6\) smaller than the tensor
product grid \(\Pattern(\mat{M}_c)\) having the same number of points.

In total, such a construction is also possible for other integer values, though the factorization
might not be that easy to find. Other normal vector directions can be approximated,
e.g.~applying the dyadic decomposition in frequency as in~\cite[Chapter 4]{Bergmann2013Thesis}.
Then, the direction to be approximated by the factorization of the matrix \(\mat{M}\)
has to be orthogonal to the direction, which should be sampled the most dense,
i.e.~for the laminate this direction of interest in the Fourier domain is
along the laminate.
%
%
%
%

\section{Summary and Conclusion}

This article generalizes the discretization of the Lippmann-Schwinger equation
and the resulting numerical algorithm to anisotropic sampling lattices. This
allows to refine other directions than the coordinate axes, even supporting
non-orthogonal refinement. This leads to smaller errors in the strain field and a better approximation of the effective stiffness tensor when taking the anisotropic properties of a material into account. Furthermore,
the orientation of the sine functions for the real-valued
discrete Fourier transform can be chosen. This allows for alignment of the ansatz functions with interfaces and increases their resolution while reducing
the Gibbs phenomenon. Especially regions and directions of high strain can thus
be resolved better. We show that these additional choices can not be
reproduced by linear transformations of the problem, e.g.~rotations of the
geometry.
Subsampling on suitable patterns makes these techniques also accessible for
data given on standard Cartesian grids present in many applications.

The application of the corresponding fast Fourier transform on patterns does
not increase the computational effort and might even be computationally
advantageous in case of lattices of rank 1. In this case the Fast Fourier
transform reduces to a one-dimensional transform, greatly reducing the computational
complexity required of the FFT algorithm. 
Other modifications applied to the Basic scheme of Moulinec and Suquet can still
be applied to the anisotropic lattice version introduced in this paper. These
modifications include adaptions of the numerical algorithm to increase both
robustness and speed. Furthermore, schemes stemming from finite difference or
finite element methods can also be incorporated into the anisotropic setting.

An open problem that has to be studied in detail is an automatically performed choice of the pattern matrix \(\mat{M}\). While the currently chosen matrices already stem from geometric interpretation of how to choose the directions of interest in the sampling lattice, the choice is up to now manually done. An analysis of the main directions of interfaces may provide a good selection of a pattern, but there may be additional restraints to take into account when selecting a sampling scheme.

\paragraph{Acknowledgment} The authors would like to thank Bernd Simeon and
Gabriele Steidl for their idea to start this collaboration.
{\small\setlength{\parskip}{.5\baselineskip}
\bibliographystyle{model1b-num-names-mod}
\bibliography{references}

\begin{thebibliography}{31}
\expandafter\ifx\csname natexlab\endcsname\relax\def\natexlab#1{#1}\fi
\providecommand{\url}[1]{\texttt{#1}}
\providecommand{\href}[2]{#2}
\providecommand{\path}[1]{#1}
\providecommand{\DOIprefix}{doi:}
\providecommand{\ArXivprefix}{arXiv:}
\providecommand{\URLprefix}{URL: }
\providecommand{\Pubmedprefix}{pmid:}
\providecommand{\doi}[1]{\href{http://dx.doi.org/#1}{\path{#1}}}
\providecommand{\Pubmed}[1]{\href{pmid:#1}{\path{#1}}}
\providecommand{\bibinfo}[2]{#2}
\ifx\xfnm\relax \def\xfnm[#1]{\unskip,\space#1}\fi
\bibitem[{{\AA}hlander and Munthe-Kaas(2005)}]{AhlanderMunthe-Kaas2005}
\bibinfo{author}{K.~{\AA}hlander}, \bibinfo{author}{H.~Munthe-Kaas},
  \bibinfo{title}{Applications of the generalized fourier transform in
  numerical linear algebra}, \bibinfo{journal}{BIT Num. Math.}
  \bibinfo{volume}{45} (\bibinfo{year}{2005}) \bibinfo{pages}{819--850}.
  \DOIprefix\doi{10.1007/s10543-005-0030-3}.
\bibitem[{Bergmann(2013{\natexlab{a}})}]{Bergmann2013FFT}
\bibinfo{author}{R.~Bergmann}, \bibinfo{title}{The fast {F}ourier transform and
  fast wavelet transform for patterns on the torus}, \bibinfo{journal}{Appl.
  Comp. Harmon. Anal.} \bibinfo{volume}{35}
  (\bibinfo{year}{2013}{\natexlab{a}}) \bibinfo{pages}{39--51}.
  \DOIprefix\doi{10.1016/j.acha.2012.07.007}.
\bibitem[{Bergmann(2013{\natexlab{b}})}]{Bergmann2013Thesis}
\bibinfo{author}{R.~Bergmann}, \bibinfo{title}{{Translationsinvariante
  R{\"a}ume multivariater anisotroper Funktionen auf dem Torus}}, Dissertation,
  Universit{\"a}t zu L{\"u}beck, \bibinfo{year}{2013}{\natexlab{b}}. \URLprefix
  \url{http://www.math.uni-luebeck.de/mitarbeiter/bergmann/publications/diss_bergmann.pdf}.
\bibitem[{Bergmann(2014)}]{Bergmann2013MPAWL}
\bibinfo{author}{R.~Bergmann}, \bibinfo{title}{The multivariate periodic
  anisotropic wavelet library}, \bibinfo{year}{2014}. \URLprefix
  \url{http://library.wolfram.com/infocenter/MathSource/8761/}.
\bibitem[{Bergmann and
  Prestin(2014{\natexlab{a}})}]{BergmannPrestin2014Interpolation}
\bibinfo{author}{R.~Bergmann}, \bibinfo{author}{J.~Prestin},
  \bibinfo{title}{Multivariate anisotropic interpolation on the torus}, in:
  \bibinfo{booktitle}{Approximation Theory XIV: San Antonio 2013},
  \bibinfo{publisher}{Springer International Publishing},
  \bibinfo{address}{Cham}, \bibinfo{year}{2014}{\natexlab{a}}, pp.
  \bibinfo{pages}{27--44}. \DOIprefix\doi{10.1007/978-3-319-06404-8_3}.
\bibitem[{Bergmann and Prestin(2014{\natexlab{b}})}]{BergmannPrestin2014dlVP}
\bibinfo{author}{R.~Bergmann}, \bibinfo{author}{J.~Prestin},
  \bibinfo{title}{Multivariate periodic wavelets of de la {V}all{\'e}e
  {P}oussin type}, \bibinfo{journal}{J. Fourier. Anal. Appl.}
  \bibinfo{volume}{21} (\bibinfo{year}{2014}{\natexlab{b}})
  \bibinfo{pages}{342--369}. \DOIprefix\doi{10.1007/978-3-319-06404-8_3}.
\bibitem[{de~Boor et~al.(1993)de~Boor, H{\"o}llig and
  Riemenschneider}]{deBoorHoelligRiemenscheider1993BoxSplines}
\bibinfo{author}{C.~de~Boor}, \bibinfo{author}{K.~H{\"o}llig},
  \bibinfo{author}{S.~Riemenschneider}, \bibinfo{title}{Box Splines},
  \bibinfo{publisher}{Springer-Verlag}, \bibinfo{address}{New York},
  \bibinfo{year}{1993}. \DOIprefix\doi{10.1007/978-1-4757-2244-4}.
\bibitem[{Chui and Li(1994)}]{ChuiLi:1994}
\bibinfo{author}{C.K. Chui}, \bibinfo{author}{C.~Li}, \bibinfo{title}{A general
  framework of multivariate wavelets with duals}, \bibinfo{journal}{Appl. Comp.
  Harmon. Anal.} \bibinfo{volume}{1} (\bibinfo{year}{1994})
  \bibinfo{pages}{368--390}. \DOIprefix\doi{10.1006/acha.1994.1023}.
\bibitem[{Galipeau and Casta{\~n}eda(2013{\natexlab{a}})}]{Galipeau2013b}
\bibinfo{author}{E.~Galipeau}, \bibinfo{author}{P.P. Casta{\~n}eda},
  \bibinfo{title}{A finite-strain constitutive model for magnetorheological
  elastomers: magnetic torques and fiber rotations}, \bibinfo{journal}{J. Mech.
  Phys. Solids} \bibinfo{volume}{61} (\bibinfo{year}{2013}{\natexlab{a}})
  \bibinfo{pages}{1065--1090}. \DOIprefix\doi{10.1016/j.jmps.2012.11.007}.
\bibitem[{Galipeau and Casta{\~n}eda(2013{\natexlab{b}})}]{Galipeau2013}
\bibinfo{author}{E.~Galipeau}, \bibinfo{author}{P.P. Casta{\~n}eda},
  \bibinfo{title}{Giant field-induced strains in magnetoactive elastomer
  composites}, in: \bibinfo{booktitle}{Proc. R. Soc. A}, volume
  \bibinfo{volume}{469}, \bibinfo{organization}{The Royal Society}, p.
  \bibinfo{pages}{20130385}. \DOIprefix\doi{10.1098/rspa.2013.0385}.
\bibitem[{Hashin and Shtrikman(1962)}]{HashinShtrikman}
\bibinfo{author}{Z.~Hashin}, \bibinfo{author}{S.~Shtrikman}, \bibinfo{title}{On
  some variational principles in anisotropic and nonhomogeneous elasticity},
  \bibinfo{journal}{J. Mech. Phys. Solids} \bibinfo{volume}{10}
  (\bibinfo{year}{1962}) \bibinfo{pages}{335--342}.
  \DOIprefix\doi{10.1016/0022-5096(62)90004-2}.
\bibitem[{Kabel et~al.(2015)Kabel, Merkert and
  Schneider}]{KMS:2015Homogenization}
\bibinfo{author}{M.~Kabel}, \bibinfo{author}{D.~Merkert},
  \bibinfo{author}{M.~Schneider}, \bibinfo{title}{Use of composite voxels in
  {FFT}-based homogenization}, \bibinfo{journal}{Comput. Method. Appl. M.}
  \bibinfo{volume}{294} (\bibinfo{year}{2015}) \bibinfo{pages}{168--188}.
  \DOIprefix\doi{10.1016/j.cma.2015.06.003}.
\bibitem[{K{\"a}mmerer et~al.(2015{\natexlab{a}})K{\"a}mmerer, Potts and
  Volkmer}]{KaemmererPottsVolkmer2015b}
\bibinfo{author}{L.~K{\"a}mmerer}, \bibinfo{author}{D.~Potts},
  \bibinfo{author}{T.~Volkmer}, \bibinfo{title}{Approximation of multivariate
  periodic functions by trigonometric polynomials based on rank-1 lattice
  sampling}, \bibinfo{journal}{J. Complexity} \bibinfo{volume}{31}
  (\bibinfo{year}{2015}{\natexlab{a}}) \bibinfo{pages}{543--576}.
  \DOIprefix\doi{10.1016/j.jco.2015.02.004}.
\bibitem[{K{\"a}mmerer et~al.(2015{\natexlab{b}})K{\"a}mmerer, Potts and
  Volkmer}]{KaemmererPottsVolkmer2015a}
\bibinfo{author}{L.~K{\"a}mmerer}, \bibinfo{author}{D.~Potts},
  \bibinfo{author}{T.~Volkmer}, \bibinfo{title}{Approximation of multivariate
  periodic functions by trigonometric polynomials based on sampling along
  rank-1 lattice with generating vector of {K}orobov form},
  \bibinfo{journal}{J. Complexity} \bibinfo{volume}{31}
  (\bibinfo{year}{2015}{\natexlab{b}}) \bibinfo{pages}{424--456}.
  \DOIprefix\doi{10.1016/j.jco.2014.09.001}.
\bibitem[{Lahellec et~al.(2003)Lahellec, Michel, Moulinec and
  Suquet}]{Lahellec2003}
\bibinfo{author}{N.~Lahellec}, \bibinfo{author}{J.C. Michel},
  \bibinfo{author}{H.~Moulinec}, \bibinfo{author}{P.~Suquet},
  \bibinfo{title}{Analysis of inhomogeneous materials at large strains using
  fast {F}ourier transforms}, in: \bibinfo{booktitle}{IUTAM symposium on
  computational mechanics of solid materials at large strains},
  \bibinfo{organization}{Springer}, pp. \bibinfo{pages}{247--258}.
  \DOIprefix\doi{10.1007/978-94-017-0297-3_22}.
\bibitem[{Langemann and Prestin(2010)}]{LangemannPrestin2010WaveletAnalysis}
\bibinfo{author}{D.~Langemann}, \bibinfo{author}{J.~Prestin},
  \bibinfo{title}{Multivariate periodic wavelet analysis},
  \bibinfo{journal}{Appl. Comp. Harmon. Anal.} \bibinfo{volume}{28}
  (\bibinfo{year}{2010}) \bibinfo{pages}{46--66}.
  \DOIprefix\doi{10.1016/j.acha.2009.07.001}.
\bibitem[{Liebscher(2014)}]{Liebscher2014}
\bibinfo{author}{A.~Liebscher}, \bibinfo{title}{Stochastic Modelling of Foams},
  Dissertation, TU Kaiserslautern, \bibinfo{year}{2014}. \URLprefix
  \url{http://www.verlag.fraunhofer.de/bookshop/buch/Stochastic-Modelling-of-Foams/242168}.
\bibitem[{Michel et~al.(2000)Michel, Moulinec and Suquet}]{Michel2000}
\bibinfo{author}{J.~Michel}, \bibinfo{author}{H.~Moulinec},
  \bibinfo{author}{P.~Suquet}, \bibinfo{title}{A computational method based on
  augmented {L}agrangians and fast {F}ourier transforms for composites with
  high contrast}, \bibinfo{journal}{Comput. Model. Eng. Sci.}
  \bibinfo{volume}{1} (\bibinfo{year}{2000}) \bibinfo{pages}{79--88}.
  \DOIprefix\doi{10.3970/cmes.2000.001.239}.
\bibitem[{Milton(2002)}]{Milton2002}
\bibinfo{author}{G.W. Milton}, \bibinfo{title}{The theory of composites},
  volume~\bibinfo{volume}{6}, \bibinfo{publisher}{Cambridge University Press},
  \bibinfo{year}{2002}. \DOIprefix\doi{10.1017/CBO9780511613357}.
\bibitem[{Moulinec and Suquet(1994)}]{MoulinecSuquet1994}
\bibinfo{author}{H.~Moulinec}, \bibinfo{author}{P.~Suquet}, \bibinfo{title}{A
  fast numerical method for computing the linear and nonlinear mechanical
  properties of composites}, \bibinfo{journal}{C. R. Acad. Sci. II B}
  \bibinfo{volume}{318} (\bibinfo{year}{1994}) \bibinfo{pages}{1417--1423}.
\bibitem[{Moulinec and Suquet(1998)}]{MoulinecSuquet1998}
\bibinfo{author}{H.~Moulinec}, \bibinfo{author}{P.~Suquet}, \bibinfo{title}{A
  numerical method for computing the overall response of nonlinear composites
  with complex microstructure}, \bibinfo{journal}{Comput. Method. Appl. M.}
  \bibinfo{volume}{157} (\bibinfo{year}{1998}) \bibinfo{pages}{69--94}.
  \DOIprefix\doi{10.1016/s0045-7825(97)00218-1}.
\bibitem[{Potts and Volkmer(2016)}]{PottsVolkmer2015}
\bibinfo{author}{D.~Potts}, \bibinfo{author}{T.~Volkmer},
  \bibinfo{title}{Sparse high-dimensional {FFT} based on rank-1 lattice
  sampling}, \bibinfo{journal}{Appl. Comput. Harm. Anal.}
  (\bibinfo{year}{2016}). \DOIprefix\doi{10.1016/j.acha.2015.05.002},
  \bibinfo{note}{to appear}.
\bibitem[{Schneider(2015)}]{Schneider2015Convergence}
\bibinfo{author}{M.~Schneider}, \bibinfo{title}{Convergence of {FFT}-based
  homogenization for strongly heterogeneous media}, \bibinfo{journal}{Math.
  Methods Appl. Sci.} \bibinfo{volume}{38} (\bibinfo{year}{2015})
  \bibinfo{pages}{2761--2778}. \DOIprefix\doi{10.1002/mma.3259}.
\bibitem[{Schneider et~al.(2016)Schneider, Merkert and Kabel}]{Schneider2016}
\bibinfo{author}{M.~Schneider}, \bibinfo{author}{D.~Merkert},
  \bibinfo{author}{M.~Kabel}, \bibinfo{title}{{FFT}-based homogenization for
  microstructures discretized by linear hexahedral elements},
  \bibinfo{journal}{submitted}  (\bibinfo{year}{2016}).
\bibitem[{Schneider et~al.(2015)Schneider, Ospald and Kabel}]{Schneider2015}
\bibinfo{author}{M.~Schneider}, \bibinfo{author}{F.~Ospald},
  \bibinfo{author}{M.~Kabel}, \bibinfo{title}{Computational homogenization of
  elasticity on a staggered grid}, \bibinfo{journal}{Int. J. Numer. Meth. Eng.}
   (\bibinfo{year}{2015}). \DOIprefix\doi{10.1002/nme.5008}.
\bibitem[{Tome and Lebensohn(1993)}]{Lebensohn}
\bibinfo{author}{R.~Tome}, \bibinfo{author}{C.~Lebensohn}, \bibinfo{title}{A
  selfconsistent approach for the simulation of plastic deformation and texture
  development of polycrystals: application to zirconium alloys},
  \bibinfo{journal}{Acta Metall. Mater.} \bibinfo{volume}{41}
  (\bibinfo{year}{1993}) \bibinfo{pages}{2611--24}.
  \DOIprefix\doi{10.1016/0956-7151(93)90130-K}.
\bibitem[{Tran et~al.(2012)Tran, Monchiet and Bonnet}]{Tran2012}
\bibinfo{author}{T.H. Tran}, \bibinfo{author}{V.~Monchiet},
  \bibinfo{author}{G.~Bonnet}, \bibinfo{title}{A micromechanics-based approach
  for the derivation of constitutive elastic coefficients of strain-gradient
  media}, \bibinfo{journal}{Int. J. Solids. Struct.} \bibinfo{volume}{49}
  (\bibinfo{year}{2012}) \bibinfo{pages}{783--792}.
  \DOIprefix\doi{10.1016/j.ijsolstr.2011.11.017}.
\bibitem[{Van~Rietbergen et~al.(1996)Van~Rietbergen, Odgaard, Kabel and
  Huiskes}]{Rietbergen1996}
\bibinfo{author}{B.~Van~Rietbergen}, \bibinfo{author}{A.~Odgaard},
  \bibinfo{author}{J.~Kabel}, \bibinfo{author}{R.~Huiskes},
  \bibinfo{title}{Direct mechanics assessment of elastic symmetries and
  properties of trabecular bone architecture}, \bibinfo{journal}{J. Biomech.}
  \bibinfo{volume}{29} (\bibinfo{year}{1996}) \bibinfo{pages}{1653--1657}.
  \DOIprefix\doi{10.1016/S0021-9290(96)80021-2}.
\bibitem[{Vond{\v{r}}ejc et~al.(2014)Vond{\v{r}}ejc, Zeman and
  Marek}]{Vondrejc2014}
\bibinfo{author}{J.~Vond{\v{r}}ejc}, \bibinfo{author}{J.~Zeman},
  \bibinfo{author}{I.~Marek}, \bibinfo{title}{An {FFT}-based {Galerkin} method
  for homogenization of periodic media}, \bibinfo{journal}{Comput. Math. Appl.}
  \bibinfo{volume}{68} (\bibinfo{year}{2014}) \bibinfo{pages}{156--173}.
  \DOIprefix\doi{10.1016/j.camwa.2014.05.014}.
\bibitem[{Willot(2015)}]{Willot2015}
\bibinfo{author}{F.~Willot}, \bibinfo{title}{Fourier-based schemes for
  computing the mechanical response of composites with accurate local fields},
  \bibinfo{journal}{C. R. Mecanique} \bibinfo{volume}{343}
  (\bibinfo{year}{2015}) \bibinfo{pages}{232--245}.
  \DOIprefix\doi{10.1016/j.crme.2014.12.005}.
\bibitem[{Zeman et~al.(2010)Zeman, Vond{\v{r}}ejc, Nov{\'a}k and
  Marek}]{ZVNM:2010NumHom}
\bibinfo{author}{J.~Zeman}, \bibinfo{author}{J.~Vond{\v{r}}ejc},
  \bibinfo{author}{J.~Nov{\'a}k}, \bibinfo{author}{I.~Marek},
  \bibinfo{title}{Accelerating a {FFT}-based solver for numerical
  homogenization of periodic media by conjugate gradients},
  \bibinfo{journal}{J. Comput. Phys.} \bibinfo{volume}{229}
  (\bibinfo{year}{2010}) \bibinfo{pages}{8065--8071}.
  \DOIprefix\doi{10.1016/j.jcp.2010.07.010}.

\end{thebibliography}
}
\end{document}